\documentclass[10pt,doublespace]{article}
\usepackage[french,english]{babel}
\usepackage{color}

\usepackage{amsmath}
\usepackage{amssymb}
\usepackage{graphicx}
\usepackage{enumerate}
\usepackage[latin1]{inputenc}
\usepackage{latexsym}
\usepackage{graphicx,epsfig}
\def\BibTeX{{\rm B\kern-.05em{\sc i\kern-.025em b}\kern-.08em
    T\kern-.1667em\lower.7ex\hbox{E}\kern-.125emX}}

\newcommand{\N}{\mathbb{N}}
\newcommand{\Z}{\mathbb{Z}}
\newcommand{\R}{\mathbb{R}}
\newcommand{\LL}{\mathbb{L}}
\newcommand{\E}{\mathbb{E}}
\newcommand{\I}{\mathbb{I}}

\newcommand{\Cov}{\mbox{Cov}}

\newcommand{\Loi}{\mathcal{L}}

\newcommand{\tend}{\overset{D}{\underset{N\rightarrow\infty}{\longrightarrow}}}
\newcommand{\ba}{\begin{eqnarray}}
\newcommand{\ea}{\end{eqnarray}}
\newcommand{\ban}{\begin{eqnarray*}}
\newcommand{\ean}{\end{eqnarray*}}
\newcommand{\be}{\begin{equation}}
\newcommand{\ee}{\end{equation}}

\def\limiteloiN{\renewcommand{\arraystretch}{0.5}
\begin{array}[t]{c}
\stackrel{{\Loi}}{\longrightarrow} \\
{\scriptstyle N\rightarrow\infty}
\end{array}\renewcommand{\arraystretch}{1}}

\def\limiteprobaN{\renewcommand{\arraystretch}{0.5}
\begin{array}[t]{c}
\stackrel{{\cal P}}{\longrightarrow} \\
{\scriptstyle N\rightarrow\infty}
\end{array}\renewcommand{\arraystretch}{1}}

\def\limitN{\renewcommand{\arraystretch}{0.5}
\begin{array}[t]{c}
\stackrel{}{\longrightarrow} \\
{\scriptstyle N\rightarrow\infty}
\end{array}\renewcommand{\arraystretch}{1}}

\newtheorem{thm}{Theorem}
\newtheorem{rem}{Remark}
\newtheorem{cor}{Corollary}

\newtheorem{prop}{Proposition}

\newtheorem{popy}{Property}

\date{}

\newenvironment{dem}{\ \\ {\bf Proof }}
{\Box\par\medskip\noindent}
 
\begin{document}
\title{Adaptive wavelet based estimator of the memory parameter for stationary Gaussian processes}
\author{\centerline{Jean-Marc Bardet$^a$, Hatem Bibi$^a$, Abdellatif
Jouini$^b$} \\
\small {\tt bardet@univ-paris1.fr}, \small {\tt
hatem.bibi@malix.univ-paris1.fr}, \small {\tt
Abdellatif.jouini@fst.mu.tn},\\
~\\
{\small $^a$ Samos-Matisse-CES, Universit\'e Paris 1, CNRS UMR 8174,
90 rue de Tolbiac, 75013 Paris, FRANCE.}\\
 {\small $^b$
D\'epartement de Math\'ematiques, Facult\'e des Sciences de Tunis, 1060
Tunis, TUNISIE.}}
 \maketitle

\begin{abstract}
This work is intended as a contribution  to a wavelet-based adaptive
estimator of the memory parameter in the classical semi-parametric
framework for Gaussian stationary processes.  In particular we
introduce and develop the choice of a data-driven optimal bandwidth.
Moreover, we  establish a central limit theorem for the estimator of
the memory parameter with the minimax rate of convergence (up to a
logarithm factor). The quality of the estimators are  attested by
simulations.

\end{abstract}
\section{Introduction}
Let $X=(X_t)_{t\in \Z}$ be a second-order zero-mean stationary
process and  its covariogram be defined
\begin{eqnarray*}
r(t)=\E (X_0\cdot X_t),~~~\mbox{for}~t\in \Z.
\end{eqnarray*}
Assume  the spectral density $f$ of $X$, with
$$
f(\lambda)=\frac 1 {2\pi} \cdot \sum_{k \in \Z} r(k) \cdot e^{-i k},
$$
exists and represents a continuous function on $[-\pi,0)[\cup ]0,\pi]$.
Consequently, the spectral density of $X$ should satisfy
the asymptotic property,
$$
f(\lambda) \sim C \cdot \frac 1 {\lambda^{D}}~~~\mbox{when $\lambda
\to 0 $},
$$
with $D<1$ called the "memory parameter" and $C>0$. If
$D\in (0,1)$, the process $X$ is a so-called  long-memory process,
if not  $X$ is called  a short memory process (see Doukhan {\it
et al.}, 2003, for more details).\\
 This paper deals with  two
 semi-parametric frameworks which are:
\begin{itemize}
\item {\bf Assumption A1:} $X$ is a zero mean stationary
Gaussian process with spectral density satisfying
$$
f(\lambda)= |\lambda |^{-D} \cdot f^*(\lambda)~~\mbox{for all}~~
\lambda \in [-\pi,0)[\cup ]0,\pi],
$$
with $f^*(0)>0$ and $f^* \in {\cal H}(D',C_{D'})$ where $0<D'$,
$0<C_{D'}$ and
$$
{\cal H}(D',C_{D'})=\Big \{ g:[-\pi,\pi] \to \R^+~\mbox{such
that}~|g(\lambda)-g(0)|\leq C_{D'} \cdot |\lambda|^{D'}~~\mbox{for
all}~\lambda \in [-\pi,\pi] \Big \}.
$$
\item {\bf Assumption A1':} $X$ is a zero-mean stationary
Gaussian process with spectral density satisfying
$$
f(\lambda)= |\lambda |^{-D} \cdot f^*(\lambda)~~\mbox{for all}~~
\lambda \in [-\pi,0)[\cup ]0,\pi],
$$
with $f^*(0)>0$ and $f^* \in {\cal H'}(D',C_{D'})$ where $0<D'$,
$C_{D'}>0$ and
$$
{\cal H'}(D',C_{D'})=\Big \{ g:[-\pi,\pi] \to \R^+~\mbox{such
that}~g(\lambda)=g(0)+C_{D'} \, |\lambda|^{D'}+o\big (
|\lambda|^{D'}\big)~~\mbox{when}~\lambda\to 0 \Big \}.
$$
\end{itemize}
\begin{rem}
A great number of earlier works concerning the estimation of the long
range parameter in a semi-parametric framework (see for instance
Giraitis {\em et al.}, 1997, 2000) are based on Assumption A1 or equivalent
assumption on $f$. Another expression (see Robinson,
1995, Moulines and Soulier, 2003 or Moulines {\it et al.}, 2007)
is $f(\lambda)=|1-e^{i\lambda}|^{-2d} \cdot f^*(\lambda)$ with $f^*$ a
function such that $|f^*(\lambda)-f^*(0)|\leq f^*(0) \cdot
\lambda^\beta$ and $0<\beta$). It is obvious that for $\beta\leq
2$ such an assumption corresponds to Assumption A1 with
$D'=\beta$. Moreover, following arguments developed in Giraitis
{\em et al.}, 1997, 2000, if $f^*\in {\cal H}(D',C_{D'})$ with
$D'>2$ is such that $f^*$ is $s\in \N^*$ times differentiable
around $\lambda=0$ with $f^{*(s)}$ satisfying a Lipschitzian
condition of degree $0<\ell<1$ around $0$, then $D' \leq s+\ell$.So for
 our purpose, $D'$ is  a more pertinent parameter than
$s+\ell$ (which is often used in no-parametric literature).
Finally, the Assumption A1' is a necessary condition to study
the following adaptive estimator of $D$.
\end{rem}
We have ${\cal H}'(D',C_{D'}) \subset {\cal
H}(D',C_{D'})$. Fractional Gaussian noises (with $D'=2$) and
FARIMA[p,d,q] processes (with also $D'=2$) represent the first and
well known examples of processes satisfying Assumption A1'
(and therefore Assumption A1). 
\begin{rem}
In Andrews and Sun (2004), an adaptive procedure covers a more general class 
of functions than ${\cal
H}(D',C_{D'})$, {\it i.e.} ${\cal H}_{AS}(D',C_{D'})$ defined by:
$$
{\cal H}_{AS}(D',C_{D'})=\Big \{ \begin{array}{c}g:[-\pi,\pi] \to \R^+~\mbox{such
that, as $\lambda \to 0$}\\
g(\lambda)=g(0)+\sum_{i=0}^k C'_i \lambda^{2i}+ C_{D'} \, |\lambda|^{D'}+o\big (
|\lambda|^{D'}\big)~~\mbox{with}~2k<D'\leq 2k+2 \end{array} \Big \}.
$$
Unfortunately, the adaptive wavelet based estimator defined below, as local or global log-periodogram estimators, is unable to be adapted to such a class (and therefore, when $D' > 2$, its convergence rate will be the same than if 
the spectral density is included in ${\cal H}_{AS}(2,C_{2})$, at the contrary to Andrew and Sun estimator). 
  \end{rem}

This work is to provide a wavelet-based  semi-parametric estimation
of the parameter $D$. This method has been introduced by Flandrin
(1989) and numerically developed by Abry {\em et al.} (1998, 2001)
and Veitch {\em et al.} (2003). Asymptotic results are reported in
Bardet {\em et al.} (2000) and more recently in Moulines {\em et
al.} (2007). Taking into account  these papers, two points of our
work can be highlighted~: first, a central limit theorem based on
 conditions which are weaker than those in Bardet {\em et al.}
(2000). Secondly, we define  an auto-driven estimator $\tilde D_n$
of $D$ (its definition being different than in Veitch {\em et
al.}, 2003). This results in  a central limit theorem followed by
$\tilde D_n$ and this estimator is proved rate optimal up to a logarithm factor 
(see below). Below we shall develop this point.\\
~\\
Define the usual Sobolev space $\tilde W(\beta,L)$ for $\beta>0$
and $L>0$,
\begin{eqnarray*}
\tilde W(\beta,L)=\left \{ g(\lambda)=\sum_{\ell \in \Z}g_\ell
e^{2\pi i\ell \lambda} \in\LL^2([0,1])\, /\,
\sum_{\ell\in\Z}(1+|\ell|)^{\beta}
|g_\ell|<\infty~~\mbox{and}~~\sum_{\ell\in\Z} |g_\ell|^2\leq L
\right \}.
\end{eqnarray*}
Let $\psi$ be a "mother" wavelet satisfying the following assumption:
\\ \\
{\bf  Assumption $W(\infty)$~:} {\em $\psi:~\R \mapsto \R$  with
$[0,1]$-support and such that}
\begin{enumerate}
\item {\em $\psi$ is included in the Sobolev class $\tilde
W(\infty,L)$ with $L>0$};
\item {\em
$\displaystyle \int_0^1  \psi(t)\,dt=0$ and $\psi(0)=\psi(1)=0$}.
\end{enumerate}
~\\
A consequence of the first point of this Assumption is: for all
$p>0$, $\sup_{\lambda \in \R} |\widehat \psi (\lambda) |
(1+|\lambda|)^{p}<\infty$, where $\widehat \psi (u)=\int_0^1
\psi(t)\,e^{-iut}dt$ is the Fourier transform of $\psi$. A useful
consequence of the second point is $\widehat \psi (u) \sim C \, u$
for $u \to 0$ with $|C|<\infty$ a real number not depending on
$u$.\\
~\\
The function $\psi$ is a smooth compactly supported function (the
interval $[0,1]$ is meant for better readability, but the following
results can be extended to another interval) with its $m$ first
vanishing moments. If $D'\leq 2$ and $0<D<1$ in Assumptions A1,
Assumption $W(\infty)$ can be replaced by a weaker assumption: \\
~\\
{\bf  Assumption $W(5/2)$~:} {\em $\psi:~\R \mapsto \R$  with
$[0,1]$-support and such that}
\begin{enumerate}
\item {\em $\psi$ is included in the Sobolev class $\tilde
W(5/2,L)$ with $L>0$};
\item {\em
$\displaystyle \int_0^1  \psi(t)\,dt=0$ and $\psi(0)=\psi(1)=0$}.
\end{enumerate}
\begin{rem}
The choice of a  wavelet satisfying Assumption $W(\infty)$ is
quite restricted because of the required smoothness of $\psi$. For
instance, the function $\psi(t)=(t^2-t+a)\exp(-1/t(1-t))$ and
$a\simeq 0.23087577$ satisfies Assumption $W(\infty)$. The class of
"wavelet" checking Assumption $W(5/2)$ is larger. For instance,
$\psi$ can be a dilated Daubechies "mother" wavelet of order $d$
with $d \geq 6$ to ensure the smoothness of the function $\psi$.It
is also possible to apply the following theory to "essentially"
compactly supported "mother" wavelet like the Lemari\'e-Meyer
wavelet. Note that it is not necessary to choose $\psi$ being a
"mother" wavelet associated to a multi-resolution analysis of
$\mathbb{L}^2(\R)$ as in the recent paper of Moulines {\em et al.}
(2007). The whole theory can be developed without this assumption,
in which case  the choice of $\psi$ is larger.
\end{rem}
If $Y=(Y_t)_{t\in \R}$ is a continuous-time process for $(a,b)\in
\R_+^*\times \R$, the "classical" wavelet coefficient $d(a,b)$ of
the process $Y$ for the scale $a$ and the shift $b$ is
\begin{eqnarray}\label{coeff_d}
d(a,b) =\frac 1{\sqrt a} \int_{\R}  \psi(\frac{t}{a}-b)Y_t \,dt.
\end{eqnarray}
However, this formula (\ref{coeff_d}) of a wavelet coefficient
cannot be computed from a time series. The support of $\psi$ being
$[0,1]$, let us take  the following approximation of
formula (\ref{coeff_d}) and define the wavelet coefficients of
$X=(X_t)_{t\in \Z}$ by
\begin{eqnarray}\label{def:eab}
e(a,b)&=&\frac {1}{\sqrt a} \sum _{k=1}^{a}\psi ( \frac {k} a)
X_{k+ab},
\end{eqnarray}
for $(a,b)\in \N_+^*\times \Z$. Note that this approximation is the
same as the wavelet coefficient computed from Mallat algorithm for
an orthogonal discrete wavelet basis (for instance with Daubechies
mother wavelet).
\begin{rem}
Here a continuous wavelet transform is considered. The discrete
wavelet transform where $a=2^j$, in other words numerically very
interesting (using Mallat cascade algorithm) is just a particular
case. The main point in studing a continuous transform is to offer a
larger number of "scales" for computing the data-driven optimal
bandwidth (see below).
\end{rem}
Under Assumption A1, for all $b \in \Z$, the asymptotic behavior of
the variance of $e(a,b)$ is a power law in  scale $a$ (when $a\to
\infty$). Indeed, for all $a \in \N^*$, $(e(a,b))_{b\in \Z}$ is a
Gaussian stationary process and (see Section more details in
\ref{CLT}):
\begin{eqnarray}\label{asymcoeff}
\E (e^2(a,0)) \sim K_{(\psi,D)} \cdot a^{D}~~~\mbox{when}~a\to
\infty,
\end{eqnarray}
with a constant $K_{(\psi,D)}$ such that,
\begin{eqnarray}\label{Kpsi}
K_{(\psi,\alpha)}=\int_{-\infty}^\infty |\widehat \psi (u)|^2 \cdot
|u|^{-\alpha}du>0~~~\mbox{for all $\alpha <1$},
\end{eqnarray}
where $\widehat \psi$ is the Fourier transform of $\psi$ (the
existence of $K_{(\psi,\alpha)}$ is established in Section
\ref{prvs}). Note that (\ref{asymcoeff}) is also checked without
the Gaussian hypothesis in Assumption A1 (the existence of the
second moment order of $X$ is sufficient). \\
The principle of the wavelet-based estimation of $D$ is linked to
this power law $a^D$. Indeed, let $(X_1,\ldots,X_N)$ be a sampled
path of $X$ and define $\widehat T_N(a)$ a sample variance of
$e(a,.)$ obtained from an appropriate choice of shifts $b$, {\it
i.e.}
\begin{eqnarray}\label{def:TN}
\widehat T_N(a)&=& \frac 1 {[N/a]} \sum _{k=1}^{[N/a]}e^2(a,k-1).
\end{eqnarray}
Then, when $a=a_N\to \infty$ satisfies $\lim_{N\to \infty} a_N \cdot
N^{-1/(2D'+1)}=\infty$, a central limit theorem for $\log(\widehat
T_N(a_N))$ can be proved. More precisely we get
$$
\log(\widehat T_N(a_N))=D \log (a_N) + \log (f^*(0)K_{(\psi,D)}) +
\sqrt{\frac {a_N} N} \cdot \varepsilon_N,
$$
with $\varepsilon_N \limiteloiN {\cal N}(0,\sigma_{(\psi,D)}^2)$ and
$\sigma_{(\psi,D)}^2>0$. As a consequence, using different scales
$(r_1 a_N,\ldots,r_\ell a_N))$ where $(r_1,\ldots,r_\ell)\in
(\N^*)^\ell$ with $a_N$ a "large enough" scale, a linear regression
of $(\log(\widehat T_N(r_i a_N))_i$ by $(\log(r_i a_N))_i$ provides
an estimator $\widehat D(a_N)$ which satisfies at the same time a
central limit theorem with a convergence rate $\sqrt{\frac N
{a_N}}$.
~\\
But the main problem is : how to select a large enough scale $a_N$
considering  that the smaller  $a_N$, the faster the convergence
rate of $\widehat D(a_N)$. An optimal solution would be to chose
$a_N$ larger but closer to $N^{1/(2D'+1)}$, but the parameter $D'$
is supposed to be unknown. In Veitch {\em et al.} (2003), an
automatic selection procedure is proposed using a chi-squared
goodness of fit statistic. This procedure is applied successfully on
a large number of numerical examples without any theoretical
proofs however. Our present  method is close to the latter. Roughly
speaking, the "optimal" choice of scale $(a_N)$ is based on the
"best" linear regression among all the possible linear regressions
of $\ell$ consecutive points $(a,\log (\widehat T_N(a)))$, where
$\ell$ is a fixed integer number. Formally speaking, a contrast is
minimized and the chosen scale $\tilde a_N$ satisfies:
$$
\frac {\log (\tilde a_N)}{\log N} \limiteprobaN \frac 1 {2D'+1}.
$$
Thus, the adaptive estimator $\tilde D_N$ of $D$ for this scale
$\tilde a_N$ is such that~:
$$
\sqrt{\frac N {\tilde a_N}} (\tilde D_N -D) \limiteloiN {\cal
N}(0,\sigma_{D}^2),
$$
with $\sigma_D^2>0$. Consequently, the minimax rate of convergence
$N^{D'/(1+2D')}$, up to a logarithm factor, for the estimation of
the long memory parameter $D$ in this semi-parametric setting (see
Giraitis {\em et al.}, 1997) is given by $\tilde D_N$. \\
~\\
Such a rate of convergence can  also be obtained by other adaptive
estimators (for  more details see below). However, $\tilde D_N$ has
 several "theoretic" advantages: firstly, it can be applied to all
$D <-1$ and $D'>0$ (which are very general conditions covering long
and short memory, in fact larger conditions  than those  usually
 required for adaptive log-periodogram or local Whittle
estimators) with a nearly optimal convergence rate. Secondly, $\tilde
D_N$ satisfies a central limit theorem and sharp confidence
intervals for $D$ can be computed (in such a case, the asymptotic
$\sigma_D^2$ is replaced by $\sigma_{\tilde D_N}^2$, for more
details see below). Finally, under additive assumptions on $\psi$
($\psi$ is supposed to have its first $m$ vanishing moments),
$\tilde D_N$ can also be applied to a process with a
polynomial trend of degree $\leq m-1$.\\
~\\
We then give a several simulations in order to appreciate empirical properties of the
adaptive estimator $\tilde D_N$. First, using a benchmark composed
of $5$ different "test" processes satisfying Assumption A1' (see
below), the central limit theorem satisfied by
$\tilde D_N$ is empirically checked. The empirical choice of
the parameter $\ell$ is also studied. Moreover, the robustness
of $\tilde D_N$ is successfully tested. Finally, the adaptive wavelet-based
estimator is compared with several existing adaptive estimators of the memory parameter
from generated paths of the $5$ different "test"
processes (Giraitis-Robinson-Samarov adaptive local log-periodogram,
Moulines-Soulier adaptive global log-periodogram, Robinson local
Whittle, Abry-Taqqu-Veitch data-driven wavelet based,
Bhansali-Giraitis-Kokoszka FAR estimators).
The simulations results of $\tilde D_N$ are convincing.
The convergence rate of $\tilde D_N$ is often ranges among the
best of the $5$ test processes (however the Robinson local
Whittle estimator $\widehat D_R$ provides more uniformly accurate estimations of $D$). Three other numerical advantages 
are offered by the adaptive wavelet-based
estimator (and not by $\widehat D_R$). Firstly, it is a very low consuming time estimator. 
Secondly it is a very robust estimator: it is not sensitive to possible
polynomial trends and seems to be consistent in non-Gaussian cases.
Finally, the graph of the log-log regression of sample variance of
wavelet coefficients is meaningful and may lead us to model data with more
general processes like locally fractional Gaussian noise (see Bardet and Bertrand, 2007).\\
~\\
The central limit theorem for sample variance of wavelet coefficient
is subject of section \ref{CLT}.Section \ref{adaptive} is concerned
with the automatic selection of the scale as well as the asymptotic
behavior of $\tilde D_N$. Finally simulations are given in section
\ref{simul} and proofs in section \ref{prvs}.
\section{A central limit theorem for the sample variance of wavelet coefficients}\label{CLT}
The following asymptotic behavior of the variance of wavelet
coefficients is the basis of all further developments.
The first point that explains all that follows is the
\begin{popy}\label{vard}
Under Assumption A1 and Assumption $W(\infty)$, for $a \in \N^*$,
$(e(a,b))_{b\in \Z}$ is a zero mean Gaussian stationary process and
it exists $M>0$ not depending on $a$ such that, for all $a \in
\N^*$,
\begin{eqnarray}\label{EA1}
\Big |\E (e^{2}(a,0))-f^*(0) K_{(\psi,D)} \cdot a^{D}\Big | \leq M
\cdot  a^{D-D'}.
\end{eqnarray}
\end{popy}
Please see Section \ref{prvs} for the proofs.
The paper of Moulines {\em et al.} (2007) gives similar results for
multi-resolution wavelet analysis. The special case of long memory
process can also be studied with weaker Assumption $W(5/2)$,
\begin{popy}\label{vard'}
Under Assumption $W(5/2)$ and Assumption A1 with $0<D<1$ and
$0<D'\leq 2$, for $a\in \N^*$, $(e(a,b))_{b\in \Z}$ is a zero mean
Gaussian stationary process and (\ref{EA1}) holds.
\end{popy}
Two corollaries can be added to both those properties. First, under
Assumption A1' a more precise result can be established.
\begin{cor}\label{cor1}
Under:
\begin{itemize}
\item Assumption A1' and Assumption $W(\infty)$;
\item {\bf or} Assumption A1' with $0<D<1$, $0<D'\leq
2$ and Assumption $W(5/2)$;
\end{itemize}
then $(e(a,b))_{b\in \Z}$ is a zero mean Gaussian stationary process
and
\begin{eqnarray}\label{equiD'}
\E (e^{2}(a,0)) = f^*(0)\Big ( K_{(\psi,D)} \cdot a^{D}+C_{D'}
K_{(\psi,D-D')} \cdot a^{D-D'}\Big )+ o\big (a^{D-D'}\big
)~~\mbox{when}~~a\to \infty.
\end{eqnarray}
\end{cor}
This corollary is key point for
the estimation of an appropriated sequence of scale $a=(a_N)$. Indeed, when $f^* \in {\cal H'}(D',C_{D'})$, then
 $f^* \in {\cal H}(D'',C_{D''})$ for all $D''$ satisfying
$0<D''\leq D'$. Therefore, Assumption A1' is required for obtaining the optimal choice of $a_N$, {\it i.e.} $a_N \simeq N^{1/(2D'+1)}$ (see below for more
details). The following corollary generalizes the above Properties
\ref{vard} and \ref{vard'}.
\begin{cor}\label{cor2}
Properties \ref{vard} and \ref{vard'} are also checked when
the Gaussian hypothesis of $X$ is replaced by $\E X_k^2<\infty$ for all $k \in \Z$.
\end{cor}
\begin{rem}
In this paper, the Gaussian hypothesis has been taken into account
merely to insure the convergence of the sample variance
(\ref{def:TN}) of wavelet coefficients following a central limit
theorem (see below). Such a convergence can also be obtained for
more general processes using a different proof of the central limit
theorem, for instance for linear processes (see a forthcoming work).
\end{rem}
As mentioned  in the introduction, this property allows an
estimation of $D$ from a log-log regression, as soon as a consistant
estimator of $\E (e^2(a,0))$ is provided from a sample
$(X_1,\ldots,X_N)$ of the time series $X$. Define then the
normalized wavelet coefficient such that
\begin{equation}\label{dtilde}
\tilde{e}(a,b)=\frac{e(a,b)}{\big (f^*(0)K_{(\psi,D)}  \cdot a^{D}
\big )^{1/2}}~~~\mbox{for $a\in \N^*$ and $b\in \Z$.}
\end{equation}
From property \ref{vard}, it is obvious that under Assumptions A1
it exists $M'>0$ satisfying for all $a\in \N^*$,
\begin{eqnarray*}
\Big |\E (\tilde{e}^2(a,0)) -1 \Big | \leq M' \cdot \frac 1
{a^{D'}}.
\end{eqnarray*}
To use this formula to estimate $D$ by a log-log regression, an
estimator of the variance of $e(a,0)$ should be considered (let us remember that a sample $(X_1,\ldots,X_N)$ of
is supposed to be known, but 
parameters ($D$, $D'$, $C_{D'}$) are unknown). Consider
the sample variance and the normalized sample variance of the
wavelet coefficient, for $1 \leq a<N$,
\begin{eqnarray} \label{samplevar}
\widehat T_{N}(a) = \frac{1}{[\frac{N}{a}]}\sum
_{k=1}^{[\frac{N}{a}]}e^{2}(a,k-1)~~~\mbox{and}~~~
\tilde{T}_{N}(a)=\frac{1}{[\frac{N}{a}]}\sum_{k=1}^{[\frac{N}{a}]}\tilde{e}^{2}(a,k-1).
\end{eqnarray}
The following proposition specifies a central limit theorem
satisfied by $\log \tilde{T}_{N}(a)$, which provides the first step
for obtaining the asymptotic properties of the estimator by log-log
regression. More generally, the following multidimensional central
limit theorem for a vector $(\log \tilde{T}_{N}(a_i))_i$ can be
established.
\begin{prop}\label{tlclog}
Define $\ell \in \N \setminus \{0,1\}$ and $(r_1,\cdots,r_\ell)
\in (\N^*)^\ell$. Let $(a_n)_{n \in \N}$ be such that $N/a_N \limitN
\infty$ and $a_N \cdot N^{-1/(1+2D')}\limitN \infty$. Under
Assumption A1 and Assumption $W(\infty)$,
\begin{equation}\label{CLTSN}
\sqrt{\frac{N}{a_N}} \Big(\log \tilde{T}_{N}(r_ia_N)\Big)_{1\leq
i\leq \ell}\tend\mathcal{N}_{\ell}\big (0\, ;\,
\Gamma(r_1,\cdots,r_\ell,\psi,D)\big ),
\end{equation}
with $\Gamma(r_1,\cdots,r_\ell,\psi,D)=(\gamma_{ij})_{1\leq i,j\leq
\ell}$ the covariance matrix such that
\begin{eqnarray}\label{cov}
\gamma_{ij}=\frac {8(r_ir_j)^{2-D}} {K^2_{(\psi,D)}d_{ij}}
\sum_{m=-\infty}^{\infty}\Big ( \int _0 ^\infty \frac {\widehat \psi
(ur_i)\overline{\widehat \psi} (ur_j) }{u^D}\cos (u \,d_{ij} m)\, du
\Big) ^2.
\end{eqnarray}
\end{prop}
The same result under weaker assumptions on $\psi$ can be also
established when $X$ is a long memory process.
\begin{prop}\label{tlclog2}
Define $\ell \in \N \setminus \{0,1\}$ and $(r_1,\cdots,r_\ell) \in
(\N^*)^\ell$. Let $(a_n)_{n \in \N}$ be such that $N/a_N \limitN
\infty$ and $a_N \cdot N^{-1/(1+2D')}\limitN \infty$. Under Assumption
$W(5/2)$ and Assumption A1 with $D\in (0,1)$ and $D'\in(0,2)$, the
CLT (\ref{CLTSN}) holds.
\end{prop}
These results can be easily generalized for processes with
polynomial trends if $\psi$ is considered
 having its first $m$ vanishing moments. i.e,
\begin{cor}\label{cor3}
Given  the same hypothesis as in Proposition \ref{tlclog} or
\ref{tlclog2} and if $\psi$ is such that  $m \in \N\setminus
\{0,1\}$ is satisfying, $\displaystyle{ \int
t^p\psi(t)\,dt=0~~\mbox{for all}~~p\in \{0,1,\ldots,m-1\}}$ the CLT
(\ref{CLTSN}) also holds for any process $X'=(X'_t)_{t\in \Z}$ such
that for all $t\in Z$, $\E X'_t=P_m(t)$ with
$P_m(t)=a_0+a_1t+\cdots+a_{m-1}t^{m-1}$ is a polynomial function and
$(a_i)_{0\leq i\leq m-1}$ are real numbers.
\end{cor}
\section{Adaptive estimator of memory parameter
using data driven optimal scales} \label{adaptive}
The CLT (\ref{CLTSN}) implies the following CLT for the vector
$(\log {\widehat T}_{N}(r_ia_N))_i$,
\begin{equation*}\label{CLTSN2}
\sqrt{\frac{N}{a_N}} \Big(\log {\widehat T}_{N}(r_ia_N) -D\log
(r_ia_N) - \log (f^*(0)K_{(\psi,D)})  \Big)_{1\leq i\leq
\ell}\tend\mathcal{N}_{\ell}\big (0\, ;\,
\Gamma(r_1,\cdots,r_\ell,\psi,D)\big ).
\end{equation*}
and therefore,
$$
\big (\log {\widehat T}_{N}(r_ia_N)\big )_{1\leq i\leq \ell} =A_N
\cdot \left
( \begin{array}{c} D \\
K \end{array}\right ) + \frac 1 {\sqrt{N/a_N} }\, \big
(\varepsilon _i\big )_{1\leq i\leq \ell},
$$
with $\displaystyle{A_N=\left
( \begin{array}{cc}\log (r_1a_N) & 1  \\
: & : \\ \log (r_\ell a_N) & 1  \end{array}\right )}$, $K=-\log
(f^*(0) \cdot K_{(\psi,D)})$ and $\displaystyle{ \big (\varepsilon
_i\big )_{1\leq i\leq \ell}\tend\mathcal{N}_{\ell}\big (0\, ;\,
\Gamma(r_1,\cdots,r_\ell,\psi,D)\big )}$. Therefore, a log-log
regression of $\big (\widehat T_{N}(r_ia_N)\big )_{1\leq i\leq
\ell}$ on scales $\big (r_ia_N\big )_{1\leq i\leq \ell}$ provides an
estimator $\displaystyle{\Big ( \begin{array}{c} \widehat D(a_N) \\
\widehat K(a_N) \end{array}\Big )}$ of $\displaystyle{\Big  (
\begin{array}{c} D \\ K \end{array}\Big  )}$ such that
\begin{equation}\label{hatd}
\Big ( \begin{array}{c} \widehat D(a_N) \\
\widehat K(a_N) \end{array}\Big )=(A_N' \cdot A_N ) ^{-1}\cdot  A_N'
\cdot
Y_{a_N}^{(r_1,\ldots,r_\ell)}~~~\mbox{with}~~Y_{a_N}^{(r_1,\ldots,r_\ell)}=\big
(\log \widehat T_{N}(r_ia_N) )_{1\leq i\leq \ell},
\end{equation}
which satisfies the following CLT,
\begin{prop}\label{tlcd}
Under the Assumptions of the Proposition \ref{tlclog},
\begin{equation}\label{CLTH}
\sqrt{\frac{N}{a_N}} \Big(\Big ( \begin{array}{c} \widehat D(a_N) \\
\widehat K(a_N) \end{array}\Big  ) -\Big  ( \begin{array}{c} D \\
K
\end{array}\Big  )\Big) \tend\mathcal{N}_{2}(0\, ; \, (A' \cdot A ) ^{-1}\cdot  A'
\cdot \Gamma(r_1,\cdots,r_\ell,\psi,D) \cdot A \cdot  (A' \cdot A ) ^{-1}),
\end{equation}
with $\displaystyle{A=\left
( \begin{array}{cc}\log (r_1) & 1  \\
: & : \\ \log (r_\ell ) & 1  \end{array}\right )}$ and
$\Gamma(r_1,\cdots,r_\ell,\psi,D)$ given by (\ref{cov}).
\end{prop}
Moreover, under Assumption A1' and if $D\in (-1,1)$, $\widehat
D(a_N)$ is a semi-parametric estimator of $D$ and its asymptotic
mean square error can be minimized with an appropriate scales
sequence $(a_N)$ reaching the well-known minimax rate of convergence
for memory parameter $D$ in this semi-parametric setting (see for
instance Giraitis {\em et al.}, 1997 and 2000). Indeed,
\begin{prop}\label{hatD}
Let $X$ satisfy Assumption A1' with $D \in (-1,1)$ and $\psi$ the
assumption $W(\infty)$. Let $(a_N)$ be a sequence such that
$a_N=[N^{1/(1+2D')}]$. Then, the estimator $\widehat D(a_N)$ is rate
optimal in the minimax sense, {\em i.e.}
$$
\limsup _{N \to \infty}\sup_{D\in (-1,1)} ~~\sup_{f^*  \in {\cal
H}(D',C_{D'})}N^{\frac {2D'}{1+2D'}} \cdot \E [{\widehat D}(a_N)
-D)^2]<+\infty.
$$
\end{prop}
\begin{rem}
As far as we know, there are no theoretic results of optimality in
case of $D\leq -1$, but according to the usual following
non-parametric theory, such minimax results can also be obtained.
Moreover, in case of long-memory processes (if $D\in (0,1)$), under
Assumption A1' for $X$ and Assumption $W(5/2)$ for $\psi$, the
estimator $\widehat D(a_N)$ is also rate optimal in the minimax
sense.
\end{rem}
In the previous Propositions \ref{tlclog} and \ref{tlcd}, the rate
of convergence of scale $a_N$ obeys to the following condition,
$$
\frac N {a_N} \limitN \infty~~\mbox{and}~~\frac {a_N}{N^{1/(1+2D')}}
\limitN \infty ~~\mbox{with}~~D'\in (0,\infty).
$$
Now, for better readability, take  $a_N=N^\alpha$. Then, the above
condition goes as follow:
\begin{eqnarray}\label{alpha}
a_N=N^\alpha~~\mbox{with}~~\alpha^* < \alpha
<1~~\mbox{and}~~\alpha^*=\frac 1 {1+2D'}.
\end{eqnarray}
Thus an optimal choice (leading to a faster convergence rate of the
estimator) is obtained for $\alpha=\alpha^*+\varepsilon$ with $\varepsilon \to 0+$. 
But $\alpha^*$ depends on $D'$ which is
unknown. To solve this problem, Veitch {\it et al.} (2003) suggest a
chi-square-based test  (constructed from a distance between the
regression line and the different points $(\log \widehat
T_N(r_ia_N), \log (r_ia_N)$). It seems to be an efficient and
interesting numerical way to estimate $D$, but without theoretical
proofs (contrary to global or local log-periodogram
procedures which are proved to reach the minimax convergence rate,
see for instance Moulines and Soulier, 2003).\\
~\\
We suggest a new procedure for the data-driven selection of optimal
scales, {\it i.e.} optimal $\alpha$. Let us consider an important parameter,
the number of considered scales $\ell \in \N \setminus \{0,1,2\}$
and set $(r_1,\ldots,r_\ell)=(1,\ldots,\ell)$. For $\alpha \in
(0,1)$, define also
\begin{itemize}
\item the vector $Y_N(\alpha)=\big (\log \widehat T_N(i \cdot N^\alpha)  \big
)_{1 \leq i \leq \ell}$;
\item the matrix
$\displaystyle{A_N(\alpha)=\left
( \begin{array}{cc}\log (N^\alpha) & 1  \\
: & : \\ \log (\ell \cdot N^\alpha) & 1  \end{array}\right )}$;
\item the contrast,
$\displaystyle
Q_N(\alpha,D,K)= \Big (Y_N(\alpha)-A_N(\alpha) \cdot \Big(
\begin{array}{c} D \\ K \end{array}\Big ) \Big )' \cdot
\Big (Y_N(\alpha)-A_N(\alpha) \cdot \Big( \begin{array}{c} D \\
K \end{array}\Big ) \Big ).
$
\end{itemize}
$Q_N(\alpha,D,K)$ corresponds to a squared distance between the
$\ell$ points $\big (\log (i \cdot N^\alpha)\,,\,\log T_N(i \cdot
N^\alpha) \big )_i$ and a line. The point is to minimize this
contrast for these three parameters. It is obvious that for a fixed
$\alpha\in (0,1)$ $Q$ is minimized from the previous least square
regression and therefore,
$$
Q_N(\widehat \alpha_N,\widehat D(a_N), \widehat K(a_N))=\min_{\alpha
\in (0,1),D<1,K \in \R} Q_N(\alpha,D,K).
$$
with $(\widehat D(a_N), \widehat K(a_N))$ obtained as in relation
(\ref{hatd}). However, since $\widehat \alpha_N$ has to be obtained from
numerical computations, the interval $(0,1)$ can be discretized as
follows,
$$
\widehat \alpha_N \in {\cal A}_N=\Big \{\frac {2}{\log
N}\,,\,\frac { 3}{\log N}\,, \ldots,\frac {\log [N/\ell]}{\log N}
\Big \}.
$$
Hence, if $\alpha \in {\cal A}_N$, it exists $k \in \{2,3,\ldots,
\log [N/\ell]\}$ such that $k=\alpha \cdot \log N$.
\begin{rem}\label{defAn}
This choice of discretization is implied by the following proof of
the consistency of $\widehat \alpha_N$. If the interval $(0,1)$ is
stepped in $N^\beta$ points, with $\beta>0$, the used proof cannot
attest this consistency. Finally, it is the same framework as the
usual discrete wavelet transform  (see for instance Veitch {\it et
al.}, 2003) but less restricted since $\log N$ may be replaced in
the previous expression of ${\cal A}_N$ by any negligible function
of $N$ compared to functions $N^\beta$ with $\beta>0$ (for instance,
$(\log N)^d$ or $d\,\log N$ can be used).
\end{rem}
Consequently, take
$$
\widehat Q_N(\alpha)=Q_N(\alpha,\widehat D(a_N), \widehat K(a_N));
$$
then, minimize $Q_N$ for variables $(\alpha,D,K)$ is
equivalent to minimize $\widehat Q_N$ for variable
$\alpha \in {\cal A}_N$, that is
\begin{eqnarray*}
\widehat Q_N(\widehat \alpha_N )=\min_{\alpha \in {\cal A}_N}
\widehat Q_N(\alpha).
\end{eqnarray*}
From this central limit theorem derives
\begin{prop}\label{hatalpha}
Let $X$ satisfy Assumption A1' and $\psi$ Assumption $W(\infty)$
(or Assumption $W(5/2)$ if $0<D<1$ and $0<D'\leq 2$). Then,
$$
\widehat \alpha_N =\frac {\log {\widehat a_N}}{\log N}
\limiteprobaN \alpha^*=\frac 1 {1+2D'}.
$$
\end{prop}
This proves also the consistency of an estimator $\widehat {D'}_N$ of
the parameter $D'$,
\begin{cor}\label{hatD'}
Taking the hypothesis of Proposition \ref{hatalpha}, we have
$$
\widehat {D'}_N=\frac {1-\widehat \alpha_N}{2 \widehat \alpha_N}
 \limiteprobaN D'.
$$
\end{cor}
The estimator $\widehat \alpha_N$ defines the selected scale
$\widehat a_N$ such that $\widehat a_N=N^{\widehat \alpha_N}$. From
a straightforward application of the proof of Proposition \ref{hatalpha} (see
the details in the proof of Theorem \ref{tildeD}), the asymptotic
behavior of $\widehat a_N$ can be specified, that is,
\begin{eqnarray}\label{hathatD}
\Pr \Big ( \frac {N^{\alpha^*}} {(\log N) ^\lambda}\leq
N^{\widehat \alpha_N} \leq N^{\alpha^*}\cdot ( \log N )^\mu \Big
)\limiteprobaN 1,
\end{eqnarray}
for all positive real numbers $\lambda$ and $\mu$ such that $\lambda
> \frac 2 {(\ell-2)D'}$ and $\mu >\frac {12}{\ell-2}$. Consequently, the selected scale is asymptotically
equal to $N^{\alpha^*}$ up to a logarithm factor. \\
~\\
Finally, Proposition \ref{hatalpha} can be used to define an
adaptive estimator of $D$. First, define the straightforward
estimator $$ \widehat {\widehat D}_N=\widehat D(\widehat a_N),
$$
which should minimize the mean square error using $\widehat a_N$.
However, the estimator $\widehat {\widehat D}_N$ does not attest a
CLT since $\Pr(\widehat \alpha_N\leq
\alpha^*)>0$ and therefore it can not be asserted  that
$\E(\sqrt{N/\widehat a_N}(\widehat {\widehat D}_N-D) )= 0$. To
establish a CLT satisfied by an adaptive estimator
$\tilde D_N$ of $D$, an adaptive scale sequence $(\tilde
a_N)=(N^{\tilde \alpha_N})$ has to be defined to ensure $\Pr(\tilde
\alpha_N \leq \alpha^*)\limitN 0$. The following theorem provides
the asymptotic behavior of such an estimator,
\begin{thm}\label{tildeD}
Let $X$ satisfy Assumption A1' and $\psi$ Assumption $W(\infty)$ (or
Assumption $W(5/2)$ if $0<D<1$ and $0<D'\leq 2$). Define,
$$
\tilde \alpha_N=\widehat \alpha_N+ \frac 3 {(\ell-2)\widehat
{D'}_N } \cdot \frac {\log \log N}{\log N},~~~\tilde a_N=N^{\tilde
\alpha_N}=N^{\widehat \alpha_N}\cdot \big ( \log N \big )^{\frac 3
{(\ell-2)\widehat {D'}_N }} ~~~\mbox{and}~~~\tilde D_N =\widehat
D(\tilde a_N).
$$
\mbox{Then, with } $\displaystyle{\sigma^2_D=(1~0)\cdot  (A' \cdot A
) ^{-1}\cdot A' \cdot \Gamma(1,\cdots,\ell,\psi,D) \cdot A \cdot (A'
\cdot A ) ^{-1}\cdot (1~0)'}$,
\begin{eqnarray}\label{CLTD2}
\sqrt{\frac{N}{N^{\tilde \alpha_N }}} \big(\tilde D_N - D \big)
\tend \mathcal{N}(0\, ; \,\sigma^2_D)~~~\mbox{and}~~~\forall
\rho>\frac {2(1+3D')}{(\ell-2)D'},~\mbox{}~~ \frac {N^{\frac
{D'}{1+2D'}} }{(\log N)^\rho} \cdot \big|\tilde D_N - D \big|
\limiteprobaN  0.
\end{eqnarray}
\end{thm}
\begin{rem}
Both the adaptive estimators $\widehat {\widehat D}_N$ and $\tilde
D_N$ converge to $D$ with a rate of convergence rate equal to the
minimax rate of convergence $N^{\frac{D'}{1+2D'}}$ up to a logarithm
factor (this result being  classical within  this semi-parametric
framework). Unfortunately, our method cannot prove that the mean
square error of both these estimators reaches the optimal rate and
therefore to be oracles.
\end{rem}
To conclude this theoretic approach, the main properties satisfied
by the estimators $\widehat {\widehat D}_N$ and $\tilde D_N$ can be
summarized as follows:
\begin{enumerate}
\item Both the adaptive estimators $\widehat {\widehat D}_N$ and $\tilde
D_N$ converge at $D$ with a rate of convergence rate equal to the
minimax rate of convergence $N^{\frac{D'}{1+2D'}}$ up to a logarithm
factor for all $D <-1$ and $D'>0$ (this being very general
conditions covering long and short memory, even larger than usual
conditions required for adaptive log-periodogram or local Whittle
estimators) whith $X$ considered a Gaussian process.
\item The estimator $\tilde D_N$ satisfies the CLT (\ref{CLTD2})
and therefore sharp confidence intervals for $D$ can be computed (in
which  case, the asymptotic matrix $\Gamma(1,\ldots,\ell,\psi,D)$ is
replaced by $\Gamma(1,\ldots,\ell,\psi,\widehat D_N)$). This is not
applicable to an  adaptive log-periodogram or local Whittle
estimators.
\item The main Property \ref{vard} is also satisfied without the
Gaussian hypothesis. Therefore, adaptive estimators $\widehat
{\widehat D}_N$ and $\tilde D_N$ can also be interesting estimators
of $D$ for non-Gaussian processes like linear or more general
processes (but a CLT similar to Theorem
\ref{tlclog} has to be established...).
\item Under additive assumptions on $\psi$ ($\psi$ is supposed to have its first
$m$ vanishing moments), both  estimators $\widehat {\widehat D}_N$
and $\tilde D_N$ can also be used for a process $X$ with a
polynomial trend of degree $\leq m-1$, which again cannot be yielded
with an adaptive log-periodogram or local Whittle estimators.
\end{enumerate}
\section{Simulations}\label{simul}
The adaptive wavelet basis estimators $\widehat {\widehat D}_N$ and
$\tilde D_N$ are new estimators of the memory parameter $D$ in the
semi-parametric frame. Different estimators of this kind are also
reported in other research works to have proved optimal. In this
paper, some theoretic advantages of adaptive wavelet basis
estimators have been highlighted. But what about concrete procedure
and results of such estimators applied to an observed sample? The
following simulations will help to answer this question.\\
First, the properties (consistency, robustness, choice of the parameter
$\ell$ and mother wavelet function $\psi$) of $\widehat {\widehat
D}_N$ and $\tilde D_N$ are investigated. Secondly, in cases of
Gaussian long-memory processes (with $D \in (0,1)$ and $D'\leq 2$),
the simulation results of the estimator $\widehat {\widehat D}_N$
are compared to those obtained
with the best known semi-parametric long-memory estimators.\\
~\\
To begin with, the simulations conditions have to be specified.
The results are obtained from $100$ generated independent samples of
each process belonging to the following "benchmark". The concrete
procedures of generation of these processes are obtained from the
circulant matrix method, as detailed in Doukhan {\it et al.} (2003).
The simulations are realized for different values of $D$,
$N$ and processes which satisfy Assumption A1' and therefore
Assumption A1 (the article of Moulines {\it et al.}, 2007, gives a
lot of details on this point):
\begin{enumerate}
\item the fractional Gaussian noise (fGn) of parameter $H=(D+1)/2$ (for $-1<D<1$) and $\sigma^2=1$.
The spectral density $f_{fGn}$ of a fGn is such that $f^*_{fGn}$ is
included in ${\cal H}(2,C_{2})$ (thus $D'=2$);
\item the FARIMA[p,d,q] process with parameter $d$ such that
$d=D/2 \in (-0.5,0.5)$ (therefore $-1<D<1$), the innovation
variance $\sigma^2$ satisfying $\sigma^2=1$ and $p, q\, \in \N$.
The spectral density $f_{FARIMA}$ of such a process is such that
$f^*_{FARIMA}$ is included in the set ${\cal H}(2,C_{2})$ (thus
$D'=2$);
\item the Gaussian stationary process $X^{(D,D')}$, such that its spectral density is
\begin{eqnarray}
f_3(\lambda)=\frac 1 {\lambda^{D}}(1+\lambda^{D'})~~~\mbox{for
$\lambda \in [-\pi,\pi]$},
\end{eqnarray}
with $D \in (-\infty,1)$ and $D'\in (0,\infty)$. Therefore $f^*_3=1+\lambda^{D'}\in
{\cal H}(D',1)$ with $D'\in (0,\infty)$.
\end{enumerate}
In the long memory frame, a "benchmark" of processes is considered for $D=0.1,~0.3,~0.5,~0.7,~0.9$:
\begin{itemize}
\item fGn processes with parameters $H=(D+1)/2$ and $\sigma^2=1$;
\item FARIMA[0,d,0] processes with $d=D/2$ and standard Gaussian
innovations;
\item FARIMA[1,d,0] processes with $d=D/2$, standard Gaussian
innovations and AR coefficient $\phi=0.95$;
\item FARIMA[1,d,1] processes with $d=D/2$, standard Gaussian
innovations and AR coefficient $\phi=-0.3$ and MA coefficient
$\phi=0.7$;
\item $X^{(D,D')}$ Gaussian processes with $D'=1$.
\end{itemize}
\subsection{Properties of adaptive wavelet basis estimators from simulations}
Below, we give the different properties of the adaptive wavelet
based method. \\
~\\
{\bf Choice of the mother wavelet $\psi$:} For short memory
processes ($D\leq 0$), let the wavelet $\psi_{SM}$ be such that
$\psi_{SM}(t)=(t^2-t+a)\exp(-1/t(1-t))$ with $a\simeq 0.23087577$.
It satisfies Assumption $W(\infty)$. Lemari\'e-Meyer wavelets can be
also investigated but this will lead to
quite different theoretic studies since its support is not bounded (but "essentially" compact). \\
For long memory processes ($0<D<1$), let the mother wavelet
$\psi_{LM}$ be such that $\psi_{LM}(t)=100\cdot
t^2(t-1)^2(t^2-t+3/14)\I_{0\leq t \leq 1}$ which satisfies Assumption $W(5/2)$. Note that 
Daubechies mother wavelet or $\psi_{SM}$
lead to "similar" results (but not as good).\\
~\\
{\bf Choice of the parameter $\ell$:} This parameter is very
important to estimate the "beginning" of the linear part of the
graph drawn by points $(\log(a_i),\log \widehat T(a_i))_i$. On the
one hand, if $\ell$ is a too small a number (for instance $\ell=3$),
another small linear part of this graph (even before the "true"
beginning $N^{\alpha^*}$) may be chosen; consequently, the
$\sqrt{MSE}$ (square root of the mean square error) of $\widehat \alpha_N$ and therefore of
$\widehat{\widehat D}_N$ or $\tilde D_N$ will be too large. On the
other hand, if $\ell$ is a too large a number (for instance
$\ell=50$ for $N=1000$), the estimator $\widehat \alpha_N$ will
certainly satisfy $\widehat \alpha_N<\alpha^*$ since it will not be
possible to consider $\ell$ different scales larger than
$N^{\alpha^*}$ (if $D'=1$ therefore $\alpha'=1/3$, then
$a_N$ has to satisfy: $N/(50a_N)=20/a_N$ is a large number and $(a_N> N^{1/3}=10$; this is not really possible). Moreover, it is
possible that a "good" choice of $\ell$ depends on the "flatness"
of the spectral density $f$, {\it i.e.} on $D'$. We have proceeded to simulations for each different values of
$\ell$ (and $N$ and $D$). Only $\sqrt{MSE}$ of estimators are presented.
The results are specified in Table 1.
\\
~\\
\begin{table}[t]\label{Table1}
{\scriptsize
\begin{center}
$N=10^3 $ ~\begin{tabular}{|c|c|c|c|c|c|c|c|}
 & $\sqrt{MSE}$ & $\ell=5$ & $\ell=10$ & $\ell=15$ & $\ell=20$ & $\ell=25$ &$\displaystyle \left \{ \begin{array}{l}\ell_1=15 \\
\ell_2=\widehat \ell \end{array} \right . $  \\
\hline \hline
fGn~$(H=\frac {D+1}2)$& $\widehat {\widehat D}_N$, $\tilde D_N$ & 0.16,~0.75&0.14,~0.19&0.13,~0.17 &\bf 0.14,~0.15 &\bf 0.14,~0.15&0.15,~0.18 \\
 &  $\widehat \alpha_N$, $\tilde \alpha_N$ & 0.12,~0.32&0.07,~0.13 &0.05,~0.08&0.04,~0.05& \bf 0.04,~0.04&0.05,~0.08\\
\hline  FARIMA$(0,\frac {D}2,0)$&  $\widehat {\widehat D}_N$, $\tilde D_N$ & 0.21,~0.81&0.15,~0.20&0.14,~0.17 &\bf 0.15,~0.15 &\bf 0.15,~0.15&0.15,~0.19 \\
 &  $\widehat \alpha_N$, $\tilde \alpha_N$ & 0.14,~0.34&0.07,~0.13&0.05,~0.09&0.05,~0.06&\bf 0.04,~0.04&0.05,~0.09\\
 \hline FARIMA$(1,\frac {D}2,0)$&  $\widehat {\widehat D}_N$, $\tilde D_N$ & 0.30,~0.96&0.28,~ 0.35&\bf 0.27,~0.29 &\bf 0.29,~0.27 &0.30,~0.30& 0.31,~0.35\\
 &  $\widehat \alpha_N$, $\tilde \alpha_N$ & 0.19,~0.44&0.15,~0.24 &0.12,~0.17&0.11,~0.15&\bf 0.11,~0.12&0.12,~0.17\\
\hline FARIMA$(1,\frac {D}2,1)$&  $\widehat {\widehat D}_N$, $\tilde D_N$ & 0.60,~0.92&0.43,~0.41&0.39,~0.35 &0.36,~0.35 &0.32,~0.33 &\bf 0.21,~0.20 \\
 &  $\widehat \alpha_N$, $\tilde \alpha_N$ & 0.17,~0.38&0.11,~0.18 &0.09,~0.12&0.07,~ 0.09&\bf 0.06,~0.07&0.09,~0.12\\
\hline  $X^{(D,D')}$, $D'=1$&  $\widehat {\widehat D}_N$, $\tilde D_N$ & 0.33,~0.68&0.29,~0.28&0.27,~0.26 &0.26,~0.27 &\bf 0.25,~0.25& 0.29,~0.30 \\
 &  $\widehat \alpha_N$, $\tilde \alpha_N$ & 0.10,~0.22&\bf 0.10,~0.07 &0.11,~0.07&0.12,~0.12& 0.13,~0.13&0.11,~0.07\\
\hline
\end{tabular}
\end{center}
\begin{center}
$N=10^4$ ~\begin{tabular}{|c|c|c|c|c|c|c|c|}
 & $\sqrt{MSE}$ & $\ell=5$ & $\ell=10$ & $\ell=15$ & $\ell=20$ & $\ell=25$ &$\displaystyle \left \{ \begin{array}{l}\ell_1=15 \\
\ell_2=\widehat \ell \end{array} \right . $  \\
\hline \hline
fGn~$(H=\frac {D+1}2)$&  $\widehat {\widehat D}_N$, $\tilde D_N$ & 0.08,~0.26&0.05,~0.05&0.05,~0.05 &\bf 0.04,~0.04 &\bf 0.04,~0.04 &\bf 0.04,~0.04\\
 &  $\widehat \alpha_N$, $\tilde \alpha_N$ & 0.08,~0.22&0.05,~0.06 &\bf 0.04,~0.05&\bf 0.04,~0.05& 0.05,~0.05&\bf 0.04,~0.05\\
\hline  FARIMA$(0,\frac {D}2,0)$&  $\widehat {\widehat D}_N$, $\tilde D_N$ & 0.08,~0.31&0.06,~0.06&\bf 0.05,~0.05 &\bf 0.05,~0.05 &\bf 0.05,~0.05 &\bf 0.05,~0.05\\
 &  $\widehat \alpha_N$, $\tilde \alpha_N$ & 0.09,~0.24&0.05,~0.07 &\bf 0.04,~0.05&\bf 0.04,~0.05& 0.05,~0.05&\bf 0.04,~0.05\\
 \hline FARIMA$(1,\frac {D}2,0)$&  $\widehat {\widehat D}_N$, $\tilde D_N$ & 0.13,~0.57&0.10,~0.10&\bf 0.09,~0.08 &\bf 0.09,~0.08&0.09,~0.09& \bf 0.09,~0.08 \\
 & $\widehat \alpha_N$, $\tilde \alpha_N$ & 0.15,~0.36&0.09,~0.16 &0.08,~0.11&0.07,~ 0.09&\bf 0.06,~0.08&0.08,~0.11\\
\hline FARIMA$(1,\frac {D}2,1)$& $\widehat {\widehat D}_N$, $\tilde D_N$ & 0.22,~0.63&0.17,~0.15&0.16,~0.13 &0.15,~0.14 & 0.15,~0.14&\bf  0.09 ,~0.09\\
 &  $\widehat \alpha_N$, $\tilde \alpha_N$ & 0.16,~0.38&0.11,~0.17 &0.08,~0.11&0.07,~0.09& \bf 0.06,~0.07&0.08,~0.11\\
\hline  $X^{(D,D')}$, $D'=1$&  $\widehat {\widehat D}_N$, $\tilde D_N$ & 0.23,~0.36&0.19,~0.15&0.18,~0.17 &0.17,~0.17 &\bf 0.15,~0.14&\bf 0.15,~0.14 \\
 &  $\widehat \alpha_N$, $\tilde \alpha_N$ & 0.10,~0.18&\bf 0.12,~0.08 &0.13,~0.12&0.14,~0.14& 0.15,~0.15&0.13,~0.12\\
\hline
\end{tabular}
\end{center}
\begin{center}
$N=10^5$ ~\begin{tabular}{|c|c|c|c|c|c|c|c|}
 & $\sqrt{MSE}$& $\ell=5$ & $\ell=10$ & $\ell=15$ & $\ell=20$ & $\ell=25$&$\displaystyle \left \{ \begin{array}{l}\ell_1=15 \\
\ell_2=\widehat \ell \end{array} \right . $  \\
\hline \hline
fGn~$(H=\frac {D+1}2)$&  $\widehat {\widehat D}_N$, $\tilde D_N$ & 0.04,~0.09&0.03,~0.03&0.02,~0.03 &\bf 0.02,~0.02 &\bf 0.02,~0.02 &\bf 0.02,~0.02\\
 &  $\widehat \alpha_N$, $\tilde \alpha_N$ & 0.07,~0.16&\bf 0.06,~0.04 &0.06,~0.06&0.07,~0.07& 0.07,~0.07&0.06,~0.06\\
\hline  FARIMA$(0,\frac {D}2,0)$& $\widehat {\widehat D}_N$, $\tilde D_N$ & 0.03,~0.13&\bf 0.02,~ 0.02&\bf 0.02,~0.02 &\bf 0.02,~0.02 &\bf 0.02,~0.02&\bf 0.02,~0.02\\
 & $\widehat \alpha_N$, $\tilde \alpha_N$ & 0.07,~0.18&0.04,~0.05 &\bf 0.04,~0.03&0.04,~0.04& 0.05,~0.05&\bf 0.04,~0.03\\
 \hline FARIMA$(1,\frac {D}2,0)$& $\widehat {\widehat D}_N$, $\tilde D_N$ & 0.05,~0.25&0.05,~ 0.04& 0.04,~0.03 & 0.04,~0.03 &0.04,~0.04 &\bf 0.03, ~0.02\\
 &  $\widehat \alpha_N$, $\tilde \alpha_N$ & 0.12,~0.30&0.07,~0.12 &0.05,~0.07&0.04,~0.06& \bf 0.04,~0.05&0.05,~0.07\\
\hline FARIMA$(1,\frac {D}2,1)$& $\widehat {\widehat D}_N$, $\tilde D_N$ & 0.08,~0.30&0.06,~ 0.04& 0.05,~0.04 & 0.05,~0.04 &0.05,~0.05 &\bf 0.04, ~0.03\\
 & $\widehat \alpha_N$, $\tilde \alpha_N$ & 0.13,~0.33&0.09,~0.15 &0.08,~0.11&0.07,~0.09&\bf 0.06,~0.08&0.08,~0.11\\
\hline  $X^{(D,D')}$, $D'=1$&  $\widehat {\widehat D}_N$, $\tilde D_N$ & 0.13,~0.19&0.11,~ 0.08&0.10,~0.08 & 0.09,~0.09 &0.09,~0.09 &\bf 0.08,~0.07\\
 &  $\widehat \alpha_N$, $\tilde \alpha_N$ & 0.09,~0.15&\bf 0.10,~0.07 &0.11,~0.09&0.12,~0.11& 0.13,~0.13&0.11,~0.09\\
\hline
\end{tabular}
\end{center}}
\caption{\it Consistency of estimators $\widehat {\widehat D}_N$, $\tilde D_N$, $\widehat \alpha_N$, $\tilde \alpha_N$
following $\ell$ from simulations of the different long-memory processes of the benchmark. For each value of $N$ ($10^3$, $10^4$ and $10^5$), of $D$ ($0.1$, $0.3$, $0.5$, $0.7$ and $0.9$)  and $\ell$ ($5$, $10$, $15$, $20$, $25$ and $(15,\widehat \ell)$), $100$ independent samples of  each process are generated. The $\sqrt{MSE}$ of each estimator is obtained from a mean of $\sqrt{MSE}$ obtained for the different values of $D$.}
\end{table}
In Table 1, two phenomena can be distinguished: the detection of
$\alpha^*$ and the estimation of $D$:
\begin{itemize}
\item  To estimate $\alpha^*$, $\ell$
has to be small enough, especially because of "$D'$ close to $0$"
and so "$\alpha'$ close to $1$" is possible. However, our
simulations indicate that $\ell$ must not be too small (for
instance $\ell=5$ leads to an important MSE for $\widehat \alpha_N$
implying an important MSE for $\widehat{\widehat D}_N$) and seems to be independent 
of $N$ (cases $N=1000$ and $N=10000$ are quite similar).
Hence, our choice is {\bf $\ell_1=15$ to estimate $\alpha^*$ for any $N$}.
\item
To estimate $D$, once $\alpha^*$ is estimated, a second value
$\ell_2$ of $\ell$ can be chosen. 
We use an adaptive procedure which, roughly speaking, consists in determining the ``end'' 
of the acceptable linear zone. Firstly, we use again the same procedure than for estimating $\widehat a_N$ but 
with scales $(a_N/i)_{1\leq i\leq \ell_1}$ and $\ell_1=15$. It provides an estimator $\widehat b_N$ corresponding 
to the maximum of acceptable (for a linear regression) scales. Secondly, {\bf the adaptive number of scales $\ell_2$ is 
computed from the formula $\ell_2=\widehat \ell=[\widehat b_N/\widehat a_N]$}. 
The simulations carried out with such values of $\ell_1$ and
$\ell_2$ are detailed in Table 1.
\end{itemize}
As it may be seen in Table 1, the choice of parameters $(\ell_1=15, \ell_2=\widehat \ell)$
provides the best results for estimating $D$, almost uniformly for all processes. \\
~\\
{\bf Consistency of the estimators $\widehat \alpha_N$ and $\tilde
\alpha_N$:} the previous numerical results (here we consider
$\ell_1=15$) show that $\widehat \alpha_N$ and $\tilde \alpha_N$
converge (very slowly) to the optimal rate $\alpha^*$, that is $0.2$
for the first four processes and $1/3$ for the fifth. Figure 1 illustrates
the evolution with $N$ of the log-log plotting and
the choice of the onset of scaling.
\begin{figure}[ht]\label{Figure1}
\[
\epsfxsize 6cm \epsfysize 6cm \epsfbox{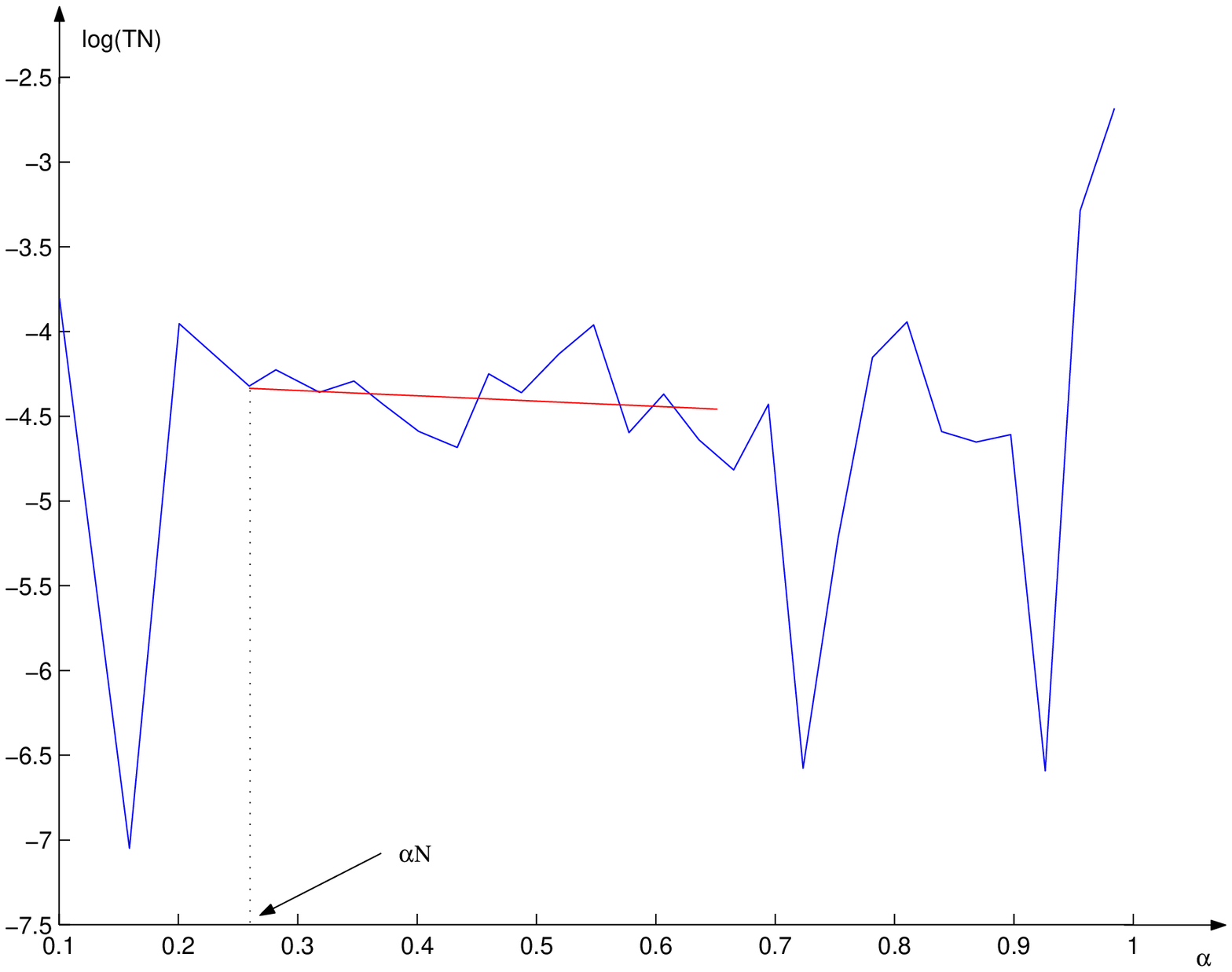}
 \hspace*{1.7
cm} \epsfxsize 6cm \epsfysize 6cm \epsfbox{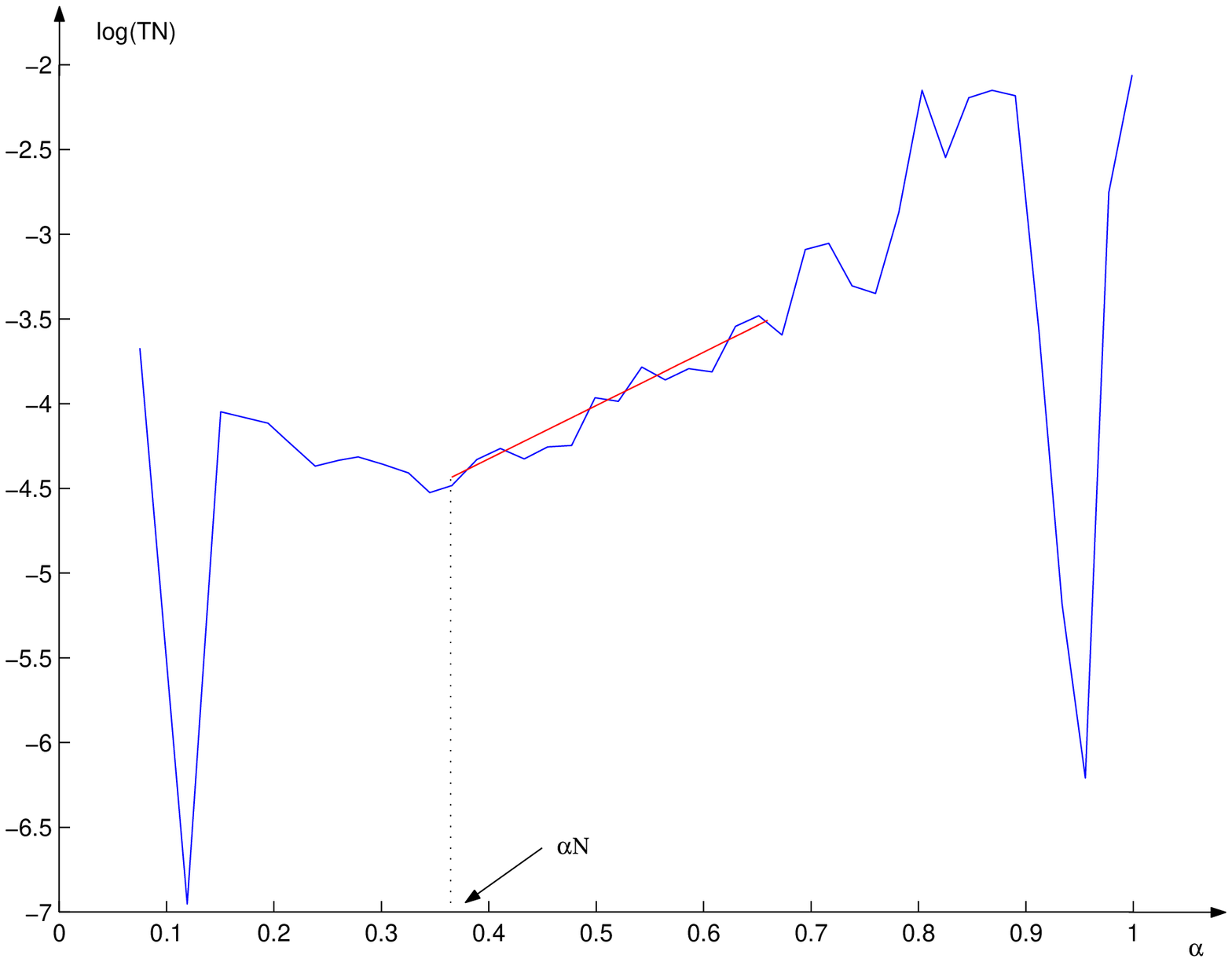}
\]
\[
\epsfxsize 6cm \epsfysize 6cm  \epsfbox{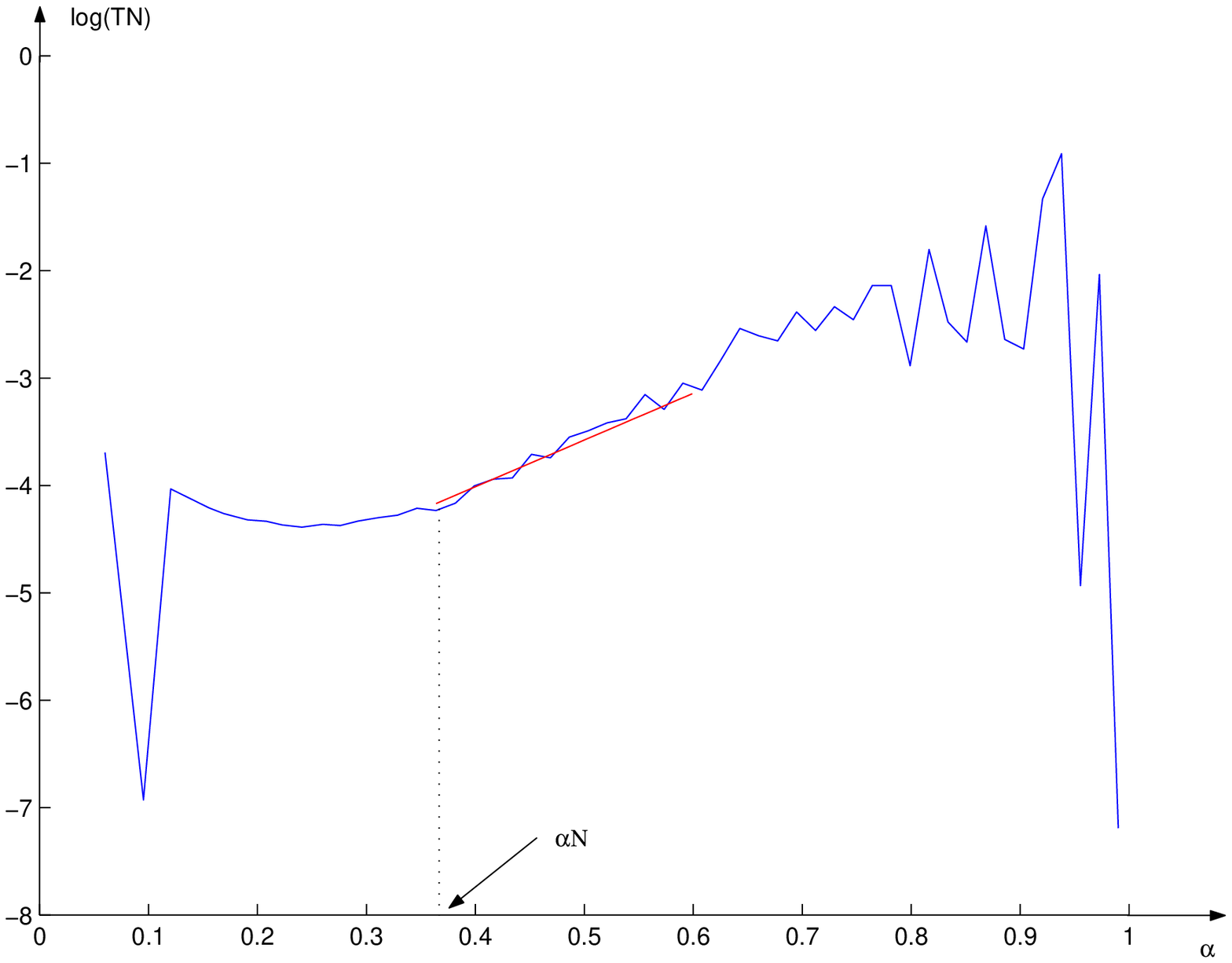}
\hspace*{1.7 cm} \epsfxsize 6cm \epsfysize 6cm
\epsfbox{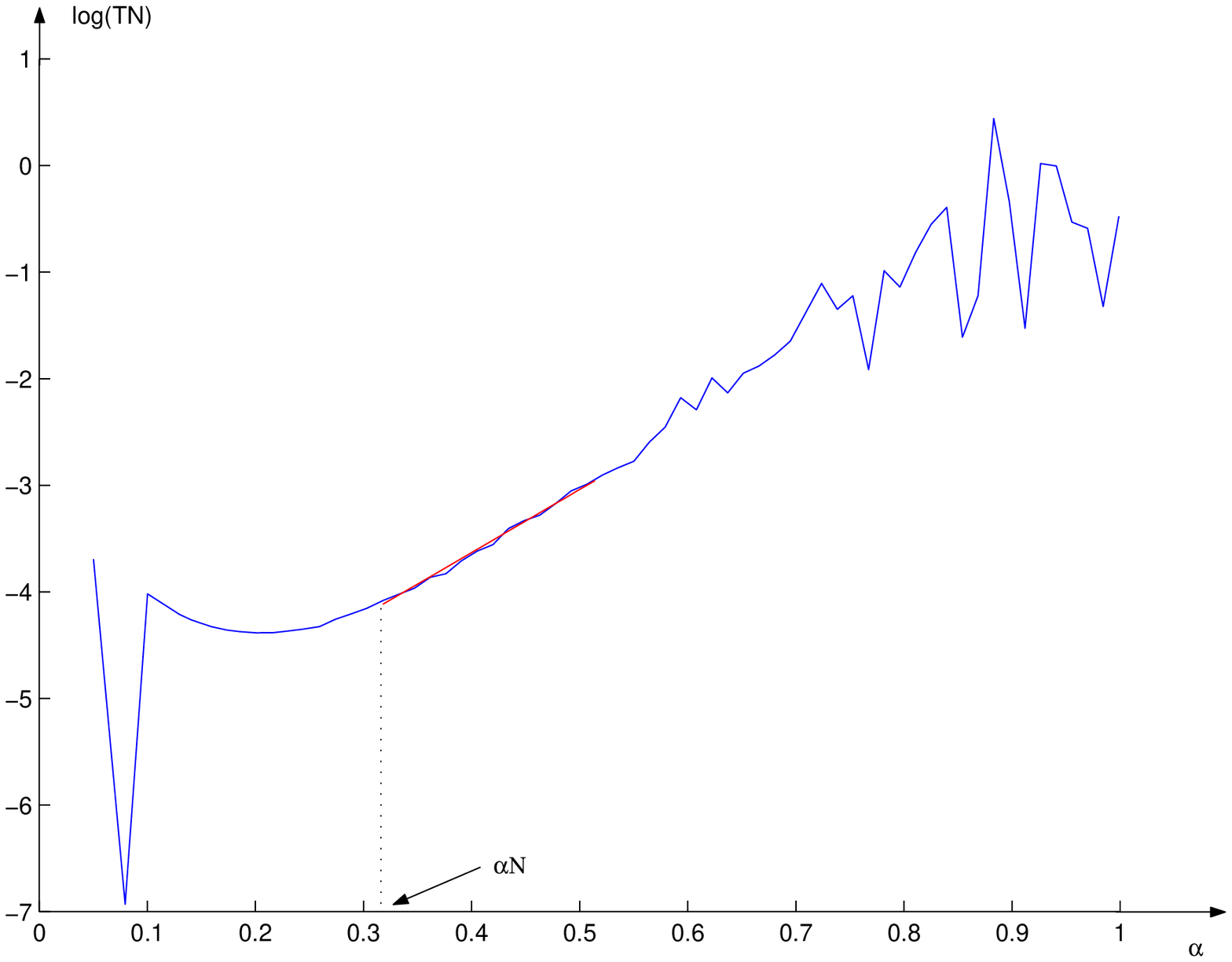}
\]
\caption{\it Log-log graphs for different samples of $X^{(D,D')}$
with $D=0.5$ and $D'=1$ when $N=10^3$ (up and left, $\widehat
{\widehat D}_N \simeq 1.04$), $N=10^4$ (up and right, $\widehat {\widehat
D}_N\simeq 0.66$), $N=10^5$ (down and left, $\widehat {\widehat
D}_N\simeq 0.62$) and $N=10^6$ (down and right,
$\widehat {\widehat D}_N\simeq 0.54$).}
\end{figure}
~\\
Figure 1 shows that $\log T_N(i \cdot N^\alpha)$ is
not a linear function of the logarithm of the scales $\log(i \cdot
N^\alpha)$ when $N$ increases and $\alpha<\alpha^*$ (a 
consequence of Property \ref{vard}: it means there is a bias).
Moreover, if $\alpha>\alpha^*$ and $\alpha$ increases, a
linear model appears with an increasing error variance.\\
~\\
{\bf Consistency and distribution of the estimators $\widehat
{\widehat D}_N$ and $\tilde D_N$:} The results of Table 1 show the
consistency with $N$ of $\widehat {\widehat D}_N$ and $\tilde D_N$
only by using $\ell_1=15$. Figure 2 provides the histograms of
$\widehat {\widehat D}_N$ and $\tilde D_N$ for $100$ independent
samples of FARIMA$(1,d,1)$ processes with $D=0.5$ and $N=10^5$.
\begin{figure}[ht]\label{Figure2}
\[
\epsfxsize 6cm \epsfysize 6cm \epsfbox{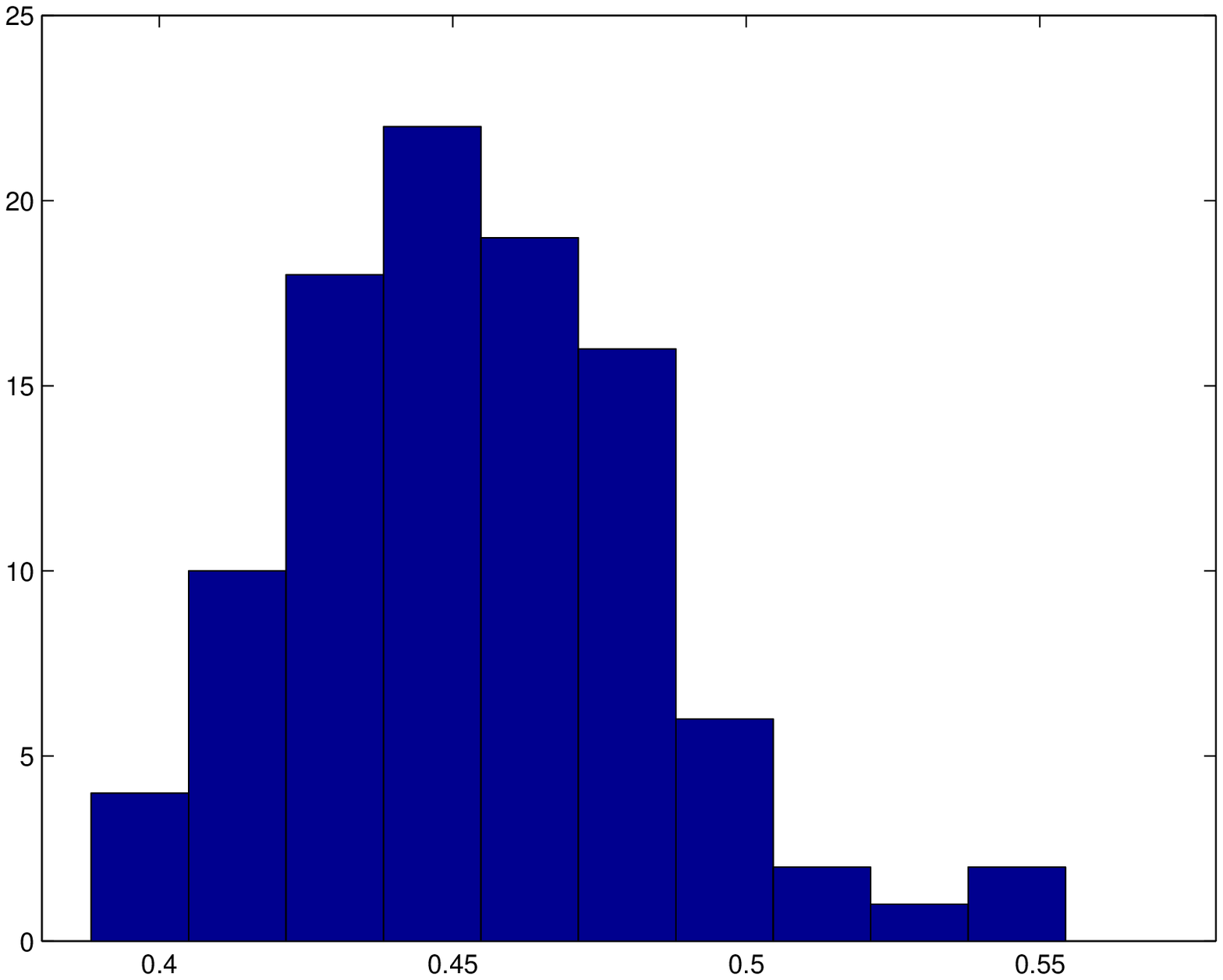}
\hspace*{1.7 cm} \epsfxsize 6cm \epsfysize 6cm
\epsfbox{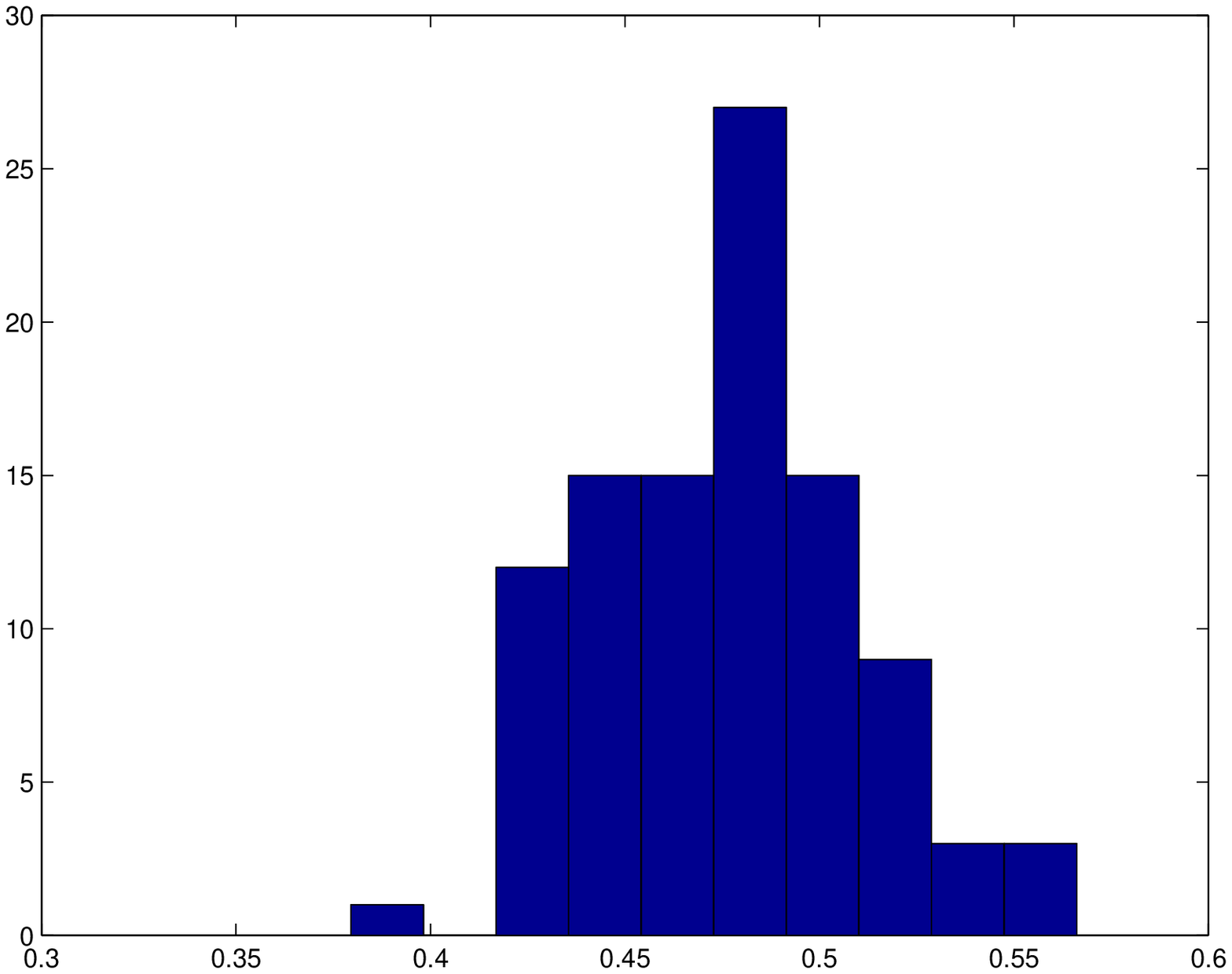}
\]
\caption{\it Histograms of $\widehat {\widehat D}_N$ and
$\tilde D_N$ for $100$ samples of FARIMA$(1,d,1)$ with $D=0.5$ for $N=10^5$.}
\end{figure}
Both the histograms of Figure 2 are similar to Gaussian distribution histograms. It is
not surprising for $\tilde D_N$ since Theorem \ref{tildeD} shows
that the asymptotic distribution of $\tilde D_N$ is a Gaussian distribution with mean
equal to $D$. The asymptotic distribution of $\widehat {\widehat D}_N$ and
the Gaussian distribution seem also to be similar. A Cramer-von Mises test of
normality indicates that both distributions of $\widehat {\widehat
D}_N$ and $\tilde D_N$ can be considered a Gaussian distribution
(respectively $W\simeq 0.07$, $p-value \simeq  0.24$
and $W \simeq  0.05$, $p-value\simeq 0.54$).\\
~\\
{\bf Consistency in case of short memory:} The following Table 2
provides the behavior of $\widehat {\widehat D}_N$ and $\tilde D_N$
if $D\leq 0$ and $D'>0$. Two processes are considered in such a
frame: a FARIMA$(0,d,0)$ process with $-0.5<d<0$ and therefore
$-1<D\leq 0$ (always with $D'=2$) and a process $X^{(D,D')}$ and $D<0$ and
$D'>0$. The results are displayed in Table \ref{Table2} (here $N=1000$, $N=10000$ and
$N=100000$, $\ell_1=15$ and $\ell_2=[5 \,N^{0.1}])$
for different choices of $D$ and $D'$. Thus it appears that
$\widehat {\widehat D}_N$ and $\tilde D_N$ can be successively
applied to short memory processes as well.
Moreover, the larger $D'$, the faster their convergence rates.\\
\begin{table}[h]\label{Table2}
{\scriptsize
\begin{center}
\begin{tabular}{|c|c|c|c|c|c|c|}
 & & FARIMA$(0,-0.25,0)$ & $X^{(-1,1)}$ & $X^{(-1,3)}$ & $X^{(-3,1)}$& $X^{(-3,3)}$\\
\hline \hline
$N=10^3$ & $\sqrt{MSE}$ $\widehat {\widehat D}_N$, $\tilde D_N$ & 0.15,~0.20&0.30,~ 0.30&0.38,~0.37 &0.36,~0.37 &0.39,~0.38\\
\hline  $N=10^4$& $\sqrt{MSE}$ $\widehat {\widehat D}_N$, $\tilde D_N$ & 0.04,~0.04&0.15,~ 0.14&0.08,~0.08 &0.13,~0.14 &0.13,~0.13\\
 \hline $N=10^5$& $\sqrt{MSE}$ $\widehat {\widehat D}_N$, $\tilde D_N$ & 0.03,~0.03&0.06,~ 0.05&0.04,~0.03 &0.04,~0.04&0.03,~0.03\\
\hline
\end{tabular}
\end{center}}
\caption {\it Estimation of the memory parameter from $100$
independent samples in case of short memory ($D\leq 0$).}
\end{table}
\newline {\bf Robustness of $\widehat {\widehat D}_N$, $\tilde D_N$:} To
conclude with the  numerical properties of the estimators, four different
processes not satisfying Assumption $A1'$ are considered:
\begin{itemize}
\item a FARIMA$(0,d,0)$ process (denoted $P1$) with innovations satisfying a uniform law (and $\E X_i^2<\infty$);
\item a FARIMA$(0,d,0)$ process (denoted $P2$) with innovations satisfying a distribution with density w.r.t.
Lebesgue measure $f(x)=3/4*(1+|x|)^{-5/2}$ for $x \in \R$ (and therefore
$\E |X_i|^2=\infty$ but $\E |X_i|<\infty$);
\item a FARIMA$(0,d,0)$ process (denoted $P3$) with innovations satisfying a Cauchy distribution (and $\E |X_i|=\infty$);
\item a Gaussian stationary process (denoted $P4$) with a spectral density $f(\lambda)=(|\lambda|-\pi/2)^{-1/2}$
for all $\lambda\in [-\pi,\pi] \setminus \{-\pi/2,\pi/2\}$. The
local behavior of $f$ in $0$ is $f(|\lambda|) \sim \sqrt{\pi/2}\,
|\lambda|^D$ with $D=0$, but the smoothness condition for $f$ in
Assumption A1 is not satisfied.
\end{itemize}
For the first $3$ processes, $D$ is varies in
$\{0.1,0.3,0.5,0.7,0.9\}$ and $100$ independent replications are
taken into account. The results of these simulations are given in Table 3.\\
\begin{table}[ht]\label{Table3}
{\scriptsize
\begin{center}
\begin{tabular}{|c|c|c|c|c|c|c|}
 & & $P1$ & $P2$ & $P3$ & $P4$\\
\hline \hline
$N=10^3$ & $\sqrt{MSE}$ $\widehat {\widehat D}_N$, $\tilde D_N$ & 0.22,~0.23&0.32,~ 0.41&0.47,~0.76&0.40,~0.41 \\
\hline  $N=10^4$& $\sqrt{MSE}$ $\widehat {\widehat D}_N$, $\tilde D_N$ & 0.06,~0.06&0.18,~ 0.28&0.24,~0.65 &0.13,~0.13 \\
\hline  $N=10^5$& $\sqrt{MSE}$ $\widehat {\widehat D}_N$, $\tilde D_N$ & 0.02,~0.02&0.02,~ 0.02&0.14,~0.47 &0.03,~0.04 \\
\hline
\end{tabular}
\end{center}}
\caption {\it Estimation of the long-memory parameter from $100$
independent samples in case of processes $P1-4$ defined above.}
\end{table}
\newline
\noindent As outlined in the theoretical part of this paper, the estimators
$\widehat {\widehat D}_N$ and $\tilde D_N$ seem also to be accurate for
$\LL^2$-linear processes. For $\LL^\alpha$-linear processes with $1\leq \alpha<2$, they
are also convergent with a slower rate of convergence. Despite the spectral density of process $P4$ does not satisfies 
the smoothness hypothesis requires in Assumptions A1 or A1', the convergence rates of $\widehat {\widehat D}_N$ and $\tilde D_N$
are still convincing. These results confirm the robustness of wavelet based estimators. 
\subsection{Comparisons with other semi-parametric long-memory parameter estimators
from simulations} 
Here we consider only long-memory Gaussian
processes ($D\in (0,1)$) based on the usual hypothesis $0<D'\leq 2$. More precisely, the "benchmark" is: $100$ generated independent
samples of each process with length $N=10^3$ and $N=10^4$ and
different values of $D$, $D=0.1,~0.3,~0.5,~0.7,~0.9$. Several
different semi-parametric estimators of $D$ are considered:
\begin{itemize}
\item $\widehat D_{BGK}$ is an "optimal" parametric Whittle estimator
obtained from a BIC criterium model selection of fractionally
differenced autoregressive models (introduced by Bhansali {it et
al.}, 2006). The required confidence interval of the estimation
$\widehat D_{BGK}$ is $[\widehat D_{R}-2/N^{1/4}\, , \, \widehat
D_{R}-2/N^{1/4}]$;
\item $\widehat D_{GRS}$ is an adaptive local periodogram estimator introduced by Giraitis
{\it et al} (2000). It requires two parameters: a bandwidth
parameter $m$, with a procedure of determination provided in this article,
 and a number of low
trimmed frequencies $l$ (satisfying different conditions but without
being fixed in this paper; after a number of  simulations,
$l=\max(m^{1/3},10)$ is chosen);
\item $\widehat D_{MS}$ is an adaptive global periodogram estimator introduced
by Moulines and Soulier (1998, 2003), also called FEXP estimator,
with bias-variance balance parameter $\kappa=2$;
\item $\widehat D_{R}$ is a local Whittle estimator introduced by
Robinson (1995). The trimming parameter is $m=N/30$;
\item $\widehat D_{ATV}$ is an adaptive wavelet based estimator
introduced by Veitch {\it et al.} (2003) using a Db4 wavelet (and
described above);
\item $\widehat {\widehat D}_N$ defined previously with $\ell_1=15$
and $\ell_2=N^{1-\widehat \alpha_N}/10$ and a
mother wavelet $\psi(t)=100\cdot t^2(t-1)^2(t^2-t+3/14)\I_{0\leq t
\leq 1}$ satisfying assumption $W(5/2)$.
\end{itemize}
Softwares (using Matlab language) for computing some of these estimators are
available on Internet (see the website of D. Veitch {\tt
http://wwww.cubinlab.ee.mu.oz.au/$\sim$darryl/} for $\widehat
D_{ATV}$ and the homepage of E. Moulines {\tt
http://www.tsi.enst.fr/$\sim$moulines/} for $\widehat D_{MS}$ and
$\widehat D_{R}$). The other softwares are available on {\tt
http://samos.univ-paris1.fr/spip/-Jean-Marc-Bardet}. Simulation results are reported in Table 4. \\
\begin{table}[p]\label{Table4} {\scriptsize
\begin{center}
$N=10^3~\longrightarrow $ ~\begin{tabular}{|c|c|c|c|c|c|c|}
 & & $D=0.1$ & $D=0.3$ & $D=0.5$ & $D=0.7$ & $D=0.9$ \\
\hline \hline
fGn~$(H=(D+1)/2)$&  $\widehat D_{BGK}$ &\bf  0.089 &0.171&0.259 &0.341&0.369 \\
&  $\widehat D_{GRS}$ & 0.114 &\bf 0.132&0.147 &0.155&0.175 \\
&  $\widehat D_{MS}$ & 0.163 &0.169&0.181 &0.195&0.191 \\
&  $\widehat D_{R}$ & 0.211 &0.220&0.215 &0.218&\bf 0.128 \\
&  $\widehat D_{ATV}$ & 0.176 &0.153&0.156 &0.164&0.162 \\
&  $\widehat {\widehat D}_N$ & 0.139 &0.147&\bf 0.133 &\bf 0.140&0.150 \\
\hline  FARIMA$(0,\frac D 2,0)$&  $\widehat D_{BGK}$ & \bf 0.094 &0.138&0.239 &0.326&0.413 \\
&  $\widehat D_{GRS}$ & 0.131&0.139&0.150 &0.150&0.162 \\
&  $\widehat D_{MS}$ & 0.172 &0.167&0.174 &0.197&0.188 \\
&  $\widehat D_{R}$ & 0.246 &0.189&0.223&0.234&0.181 \\
&  $\widehat D_{ATV}$ & 0.128 & \bf  0.107& \bf  0.081 & \bf  0.074& \bf  0.065 \\
&  $\widehat {\widehat D}_N$ & 0.161 &0.146&0.149 &0.149&0.161 \\
\hline FARIMA$(1,\frac D 2,0)$&   $\widehat D_{BGK}$ & \bf  0.146 & \bf  0.203&0.239 &0.236&0.212 \\
&  $\widehat D_{GRS}$ & 0.519 &0.545&0.588 &0.585&0.830 \\
&  $\widehat D_{MS}$ & 0.235 &0.258&0.256 &0.252&0.249\\
&  $\widehat D_{R}$ & 0.242 &0.241& \bf  0.234 & \bf  0.202 \bf  &0.144 \\
&  $\widehat D_{ATV}$ & 0.248 &0.267&0.280 &0.268&0.375 \\
&  $\widehat {\widehat D}_N$ & 0.340 &0.319&0.314 &0.315 &0.334 \\
\hline FARIMA$(1,\frac D 2,1)$&  $\widehat D_{BGK}$ & 0.204 &0.253&0.342 &0.363&0.384 \\
&  $\widehat D_{GRS}$ & 0.901 &0.894&0.866 &0.870&0.893 \\
&  $\widehat D_{MS}$ & 0.181 & \bf  0.175& \bf  0.180 & \bf  0.185&0.181 \\
&  $\widehat D_{R}$ & 0.204 &0.200&0.200 &0.191& \bf  0.130 \\
&  $\widehat D_{ATV}$ & 0.392 &0.380&0.371 &0.343&0.355 \\
&  $\widehat {\widehat D}_N$ &  \bf  0.170 &0.218&0.225 &0.226&0.213 \\
\hline  $X^{(D,D')}$, $D'=1$&  $\widehat D_{BGK}$ & \bf   0.090 & \bf  0.139&0.261 &0.328&0.388 \\
&  $\widehat D_{GRS}$ & 0.342 &0.339&0.331 &0.300&0.315 \\
&  $\widehat D_{MS}$ & 0.176 &0.178&0.182 & \bf  0.166&0.177 \\
&  $\widehat D_{R}$ & 0.219 &0.232&0.231 &0.173& \bf  0.167 \\
&  $\widehat D_{ATV}$ & 0.153 &0.161& \bf  0.168 &0.176&0.176 \\
&  $\widehat {\widehat D}_N$ & 0.284 &0.294&0.293 &0.292&0.288 \\
\hline
\end{tabular}
\end{center}}
\end{table}
\begin{table} {\scriptsize
\begin{center}
$N=10^4~\longrightarrow $ ~\begin{tabular}{|c|c|c|c|c|c|c|}
 & & $D=0.1$ & $D=0.3$ & $D=0.5$ & $D=0.7$ & $D=0.9$ \\
\hline \hline
fGn~$(H=(D+1)/2)$&  $\widehat D_{BGK}$ & 0.062 &0.143&0.182 &0.171&0.182 \\
&  $\widehat D_{GRS}$ & 0.040 &0.047&0.054 &0.068&0.066 \\
&  $\widehat D_{MS}$ & 0.069 &0.064&0.061 &0.071&0.063 \\
&  $\widehat D_{R}$ & 0.063 &0.055&0.058 &0.063&0.052 \\
&  $\widehat D_{ATV}$ & \bf  0.036 & \bf 0.042&  \bf  0.041 &0.047&0.045 \\
&  $\widehat {\widehat D}_N$ & 0.050 &0.040& \bf  0.041 & \bf  0.039& \bf  0.040 \\
\hline  FARIMA$(0,\frac D 2,0)$&  $\widehat D_{BGK}$ & 0.059 &0.141&0.195 &0.187&0.178 \\
&  $\widehat D_{GRS}$ & 0.042 &0.048&0.050 &0.046&0.057 \\
&  $\widehat D_{MS}$ & 0.072 &0.055&0.066 &0.059&0.065 \\
&  $\widehat D_{R}$ & 0.073 &0.053&0.064 &0.057&0.059 \\
&  $\widehat D_{ATV}$ &  \bf  0.026 & \bf  0.038& \bf  0.039 & \bf  0.032& \bf  0.022 \\
&  $\widehat {\widehat D}_N$ & 0.053 &0.050&0.056 &0.055&0.044 \\
\hline FARIMA$(1,\frac D 2,0)$&   $\widehat D_{BGK}$ & 0.085&0.148&0.146 &0.164&0.120 \\
&  $\widehat D_{GRS}$ & 0.179&0.175&0.182 &0.192&0.190 \\
&  $\widehat D_{MS}$ & 0.109 &0.105&0.099 &0.100&0.094 \\
&  $\widehat D_{R}$ & \bf  0.063 & \bf 0.059& \bf 0.057 & \bf 0.054& \bf 0.054 \\
&  $\widehat D_{ATV}$ & 0.118&0.101&0.088 &0.120&0.081 \\
&  $\widehat {\widehat D}_N$ & 0.095 &0.085&0.093 &0.081&0.097 \\
\hline FARIMA$(1,\frac D 2,1)$&  $\widehat D_{BGK}$ & 0.111 &0.201&0.189 &0.202&0.181 \\
&  $\widehat D_{GRS}$ & 0.308 &0.321&0.306&0.314&0.311 \\
&  $\widehat D_{MS}$ & 0.070 &0.064&0.065 & \bf 0.064&0.069 \\
&  $\widehat D_{R}$ &  \bf 0.063 & \bf 0.057& \bf 0.060 & \bf 0.064& \bf 0.052 \\
&  $\widehat D_{ATV}$ & 0.114 &0.118&0.103 &0.102&0.093 \\
&  $\widehat {\widehat D}_N$ & 0.095 &0.099&0.087 &0.101&0.090 \\
\hline  $X^{(D,D')}$, $D'=1$&  $\widehat D_{BGK}$ & 0.069 &0.110&0.204 &0.190&0.197 \\
&  $\widehat D_{GRS}$ & 0.192 &0.185&0.172 &0.177&0.190 \\
&  $\widehat D_{MS}$ & 0.083 &0.059&0.071 &0.066&0.068 \\
&  $\widehat D_{R}$ &  \bf 0.066 & \bf 0.057& \bf 0.068 & \bf 0.054& \bf 0.064 \\
&  $\widehat D_{ATV}$ & 0.124 &0.131&0.139 &0.147&0.153\\
&  $\widehat {\widehat D}_N$ & 0.158 &0.143&0.152 &0.158&0.155 \\
\hline
\end{tabular}
\end{center}}
\caption {\it Comparison of the different log-memory parameter estimators for processes of the benchmark.
For each process and value of $D$ and $N$, $\sqrt{MSE}$ are computed from $100$ independent generated samples. }
\end{table} 
\newline 
{\bf Comments on the results of Table 4:} These simulations allow to
distinguish four "clusters" of estimators.
\begin{itemize}
\item $\widehat D_{BGK}$ is obtained from a BIC-criterium hierarchical model
selection (from $2$ to $11$ parameters, corresponding to the length
of the approximation of the Fourier expansion of the spectral
density) using Whittle estimation. For these simulations, the
BIC criterion is generally minimal for $5$ to $7$ parameters to be estimated.
Simulation results are not very satisfactory except for $D=0.1$ (close to
the short memory). Moreover, this procedure is rather
time-consuming.
\item $\widehat D_{GRS}$ offers good results for fGn and
FARIMA$(0,d,0)$. However, this estimator does not converge fast
enough for the other processes.
\item Estimators $\widehat D_{MS}$ and $\widehat D_{R}$ have
similar properties. They (especially $\widehat D_{R}$) are very interesting because they offer the
same fairly good rates of convergence for all processes of the
benchmark.
\item Being built on similar principles, estimators $\widehat D_{ATV}$
and $\widehat {\widehat D}_N$ have similar behavior as well.
Their convergence rates are the fastest for fGn and FARIMA$(0,d,0)$
and are almost close to fast ones for the other processes. Their times of
computing, especially for $\widehat D_{ATV}$ for which the
computations of wavelet coefficients with that the Mallat algorithm,
are the shortest.
\end{itemize}
{\bf Conclusion:} Which estimator among those studied above has to be chosen in a practical frame, {\it i.e.} 
an observed time series?
We propose the following procedure for estimating an eventual long memory parameter: 
\begin{enumerate}
\item Firstly, since this procedure is very low time consuming and applicable to processes 
with smooth trends, draw the log-log regression of wavelet coefficients' variances onto scales. If 
a linear zone appears in this graph, consider the estimator $\widehat {\widehat D}_N$ (or $\widehat D_{ATV}$) of $D$.
\item If a linear zone appears in the previous graph and if the observed time series seems to 
be without a trend, compute $\widehat D_{R}$. 
\item Compare both the estimated value of $D$ from confidence intervals (available for $\widehat {\widehat D}_N$ or
$\widehat D_{ATV}$ and $\widehat D_{R}$).
\end{enumerate}
\section{Proofs}\label{prvs} \dem[Property \ref{vard}] The arguments of
this proof are similar to those of Abry {\it et al.} (1998) or
Moulines {\it et al.} (2007). First, for $a\in \N^*$,
\begin{eqnarray}
\nonumber \E (e^{2}(a,0))&=& \frac 1 a \sum_{k=1}^a \sum_{k'=1}^a
\psi(k/a)\psi(k'/a)\E(X_{k}X_{k'})\\
\nonumber &=&\frac 1 a \sum_{k=1}^a \sum_{k'=1}^a
\psi(k/a)\psi(k'/a)r(k-k')\\
\nonumber &=&\frac 1 a \sum_{k=1}^a \sum_{k'=1}^a
\psi(k/a)\psi(k'/a)
\int_{-\pi}^\pi f(\lambda)e^{i\lambda  (k-k')} d\lambda\\
\nonumber &=& \int_{-a\pi}^{a\pi} f\big ( \frac  u a \big ) \times
\frac 1 {a^2}\sum_{k=1}^a \sum_{k'=1}^a \psi\big ( \frac  k a \big
)\psi\big ( \frac  {k'} a \big ) e^{iu\big (\frac  k a-\frac  {k'}
a\big )}
du\\
\label{partie1} &=& \int_{-a\pi}^{a\pi} f\big ( \frac  u a \big )
\times \left \{ \Big (\frac 1 {a}\sum_{k=1}^a \psi\big ( \frac  k a
\big ) \cos \big ( \frac  k a u \big )\Big )^2 + \Big (\frac 1
{a}\sum_{k=1}^a \psi\big ( \frac  k a \big ) \sin \big ( \frac  k a
u \big )\Big )^2 \right \}\, du
\end{eqnarray}
Now, it is well known that if $\psi\in \tilde W(\beta,L)$ the
Sobolev space with parameters $\beta>1/2$ and $L>0$, then
\begin{eqnarray}\label{induction}
\sup_{|u|\leq a\pi}\Delta_a(u) \leq C_{\beta,L} \, \frac 1
{a^{\beta-1/2}}~~\mbox{with}~~\Delta_a(u):=\Big |\frac 1 a
\sum_{k=1}^a \psi\big ( \frac k a \big )e^{-iu\frac k a} -\int_0^1
\psi(t)e^{-iut}dt \Big |,
\end{eqnarray}
with $C_{\beta,L}>0$ only depending on $\beta$ and $L$ (see for
instance Devore and Lorentz, 1993). Therefore if $\psi$ satisfies
Assumption $W(\infty)$ and $X$ Assumption A1, for all $\beta>1/2$,
since $\sup_{u\in \R} |\widehat \psi (u) | <\infty$,
\begin{eqnarray}
\nonumber \Big |\E (e^{2}(a,0))-\int_{-a\pi}^{a\pi} f\big ( \frac u
a \big ) \times|\widehat \psi (u)|^2\, du \Big | &\leq &
2C_{\beta,L} \, \frac 2 {a^{\beta-3/2}} \int _{0}^{a\pi} f\big (
\frac  u a \big ) |\widehat \psi(u)|\, du+C^2_{\beta,L} \, \frac 2
{a^{2\beta-2}}\int
_{0}^{a\pi} f\big ( \frac  u a \big )\, du\\
\label{Fourier} &\leq &2 \cdot C^2_{\beta,L} \, \frac 2
{a^{2\beta-3}}\int _{0}^{\pi} f(v )\, dv,
\end{eqnarray}
since $\sup_{u\in \R} (1+u^n)|\widehat \psi (u)| <\infty$ for all $n
\in \N$. Consequently, if $\psi$ satisfies Assumption
$W(\infty)$, for all $n >0$, for all $a \in \N^*$, there exists
$C(n)>0$ not depending on $a$ such that
\begin{eqnarray}\label{Fourier2}
\Big |\E (e^{2}(a,0))-\int_{-a\pi}^{a\pi} f\big ( \frac u a \big )
\times|\widehat \psi (u)|^2\, du \Big | &\leq &C(n)\,  \frac 1
{a^{n}}.
\end{eqnarray}
But from Assumption $W(\infty)$, for all $c<1$,
$$
K_{(\psi,c)}=\int_{-\infty}^{\infty} \frac {|\widehat \psi
(u)|^2}{|u|^c} \, du <\infty,
$$
because Assumption $W(\infty)$ implies that $|\widehat \psi
(u)|=O(|u|)$ when $u \to 0$ and there exists $p>1-c$ such that $\sup
_{u \in \R}|\widehat \psi (u)|^2(1+ |u|)^{p}<\infty$. Moreover, for
all $p>1-c$,
\begin{eqnarray*}
\Big |\int_{-a\pi}^{a\pi} \frac {|\widehat \psi (u)|^2}{|u|^c} \,
du-K_{(\psi,c)}\Big |&=& 2\int_{a\pi}^{\infty} \frac {|\widehat
\psi
(u)|^2}{u^c} \, du \\
&\leq & C \cdot \int_{a\pi}^{\infty} \frac 1 {u^{p+c}} \, du \\
&\leq & C' \cdot \frac 1 {a^{p+c-1}},
\end{eqnarray*}
with $C>0$ and $C'>0$ not depending on $a$. As a consequence,
under Assumption A1, for all $p>1-D$, all $n \in \N$ and all $a
\in \N^*$,
\begin{eqnarray*}
\Big |\E (e^{2}(a,0))-f^*(0) \cdot \int_{-\infty}^{\infty} \frac
{|\widehat \psi (u)|^2}{|u/a|^{D}} \, du \Big |  & \leq & 2f^*(0)a^D
\int_{a\pi}^{\infty} \frac {|\widehat \psi (u)|^2}{u^{D}} \, du+
C_{D'}a^{D-D'} \int_{-a\pi}^{a\pi} \frac {|\widehat \psi
(u)|^2}{|u|^{D-D'} }\, du+C(n)\,  \frac 1
{a^{n}}\\
\Longrightarrow ~~\Big |\E (e^{2}(a,0))-f^*(0) K_{(\psi,D)} \cdot
a^{D}\Big |  &\leq &  C' f^*(0)\cdot a^{1-p}+C_{D'} K_{(\psi,D-D')} \cdot
a^{D-D'}.
\end{eqnarray*}
Now, by choosing $p$ such that $1-p<D-D'$, the inequality
(\ref{EA1}) is obtained. $\Box$\\
\dem[Property \ref{vard'}] Using the proof of previous Property
\ref{vard}, with Assumption $W(5/2)$, $\psi$ is included in a
Sobolev space $\tilde W(5/2,L)$, inequality (\ref{induction}) is
checked with $\beta=5/2$ and (\ref{Fourier}) is replaced by
\begin{eqnarray}\label{Fourier52}
\Big |\E (e^{2}(a,0))-a \int_{-a\pi}^{a\pi} f\big ( \frac u a \big
) \times|\widehat \psi (u)|^2\, du \Big | \leq 2 \cdot C^2_{5/2,L} \, \frac
2 {a^{2}}\int _{0}^{\pi} f(v )\, dv,
\end{eqnarray}
since $\sup_{u\in \R} (1+u^{3/2})|\widehat \psi (u)| <\infty$.
Therefore, inequality (\ref{Fourier2}) is replaced by
\begin{eqnarray*}
\Big |\E (e^{2}(a,0))-a \int_{-a\pi}^{a\pi} f\big ( \frac u a \big
) \times|\widehat \psi (u)|^2\, du \Big | &\leq &C(2)\,  \frac 1
{a^{2}}.
\end{eqnarray*}
The end of the proof is similar to the end of the previous proof,
but now $K_{(\psi,c)}$ exists for $-2< c< 1$ and
$$
\Big |\int_{-a\pi}^{a\pi} \frac {|\widehat \psi (u)|^2}{|u|^c} \,
du-K_{(\psi,c)}\Big |\leq C' \cdot \frac 1 {a^{2+c}}.
$$
Finally, under Assumption A1', for all $a \in \N^*$, since $-2<
D-D'<1$,
\begin{eqnarray*}
\Big |\E (e^{2}(a,0))-f^*(0) K_{(\psi,D)} \cdot a^{D}\Big |  &\leq &
C_{D'} K_{(\psi,D-D')} \cdot a^{D-D'}+C'\, \frac 1 {a^{2}},
\end{eqnarray*}
which achieves the proof. $\Box$ \\
\dem[Corollary \ref{cor1}] Both these proofs provide main arguments
to establish (\ref{equiD'}). For better readability , we will
consider only  Assumption A1' and Assumption $W(\infty)$ (the  long
memory process being  similar). The main difference consists in
specifying the asymptotic behavior of $\displaystyle
\int_{-a\pi}^{a\pi} f\big ( \frac u a \big ) \times|\widehat \psi
(u)|^2\, du$. But,
\begin{eqnarray}\label{som1}
\int_{-a\pi}^{a\pi} f\big ( \frac u a \big ) \times|\widehat \psi
(u)|^2\, du=\int_{-\sqrt{a}}^{\sqrt{a}} f\big ( \frac u a \big )
\times|\widehat \psi (u)|^2\, du+ 2 \int_{\sqrt{a}}^{a \pi} f\big (
\frac u a \big ) \times|\widehat \psi (u)|^2\, du.
\end{eqnarray}
The asymptotic behavior of $\widehat \psi(u)$ when $u\to \infty$
($\psi$ is considered  to satisfy Assumption $W(\infty)$), this
behavior induces that
\begin{eqnarray}\label{som2}
\int_{\sqrt{a}}^{a \pi} f\big ( \frac u a \big ) \times|\widehat
\psi (u)|^2\, du \leq C a^D \int_{\sqrt{a}}^{\infty } u^{-D}
\times|\widehat \psi (u)|^2\, du \leq \frac {C(n)} {a^n},
\end{eqnarray}
for all $n \in \N$. Moreover,
\begin{eqnarray}
\nonumber &&\int_{-\sqrt{a}}^{\sqrt{a}} f\big ( \frac u a \big )
\, |\widehat \psi (u)|^2\, du= f^*(0) \int_{-\sqrt{a}}^{\sqrt{a}}
\Big (\Big | \frac u a \Big |^{-D} +C_{D'}\Big | \frac u a \Big
|^{D'-D}  \Big ) |\widehat \psi (u)|^2\, du\\
\label{som3}&& \hspace{5cm} + \int_{-\sqrt{a}}^{\sqrt{a}}\Big (
f\big ( \frac u a \big ) - f^*(0) \Big (\Big | \frac u a \Big
|^{-D} +C_{D'}\Big | \frac u a \Big |^{D'-D}  \Big )\Big )
|\widehat \psi (u)|^2\, du.
\end{eqnarray}
From computations of previous proofs,
\begin{eqnarray}\label{som4}
\int_{-\sqrt{a}}^{\sqrt{a}} \Big (\Big | \frac u a \Big |^{-D}
+C_{D'}\Big | \frac u a \Big |^{D'-D}  \Big ) |\widehat \psi
(u)|^2\, du=K_{(\psi,D)}\cdot a^{D} +C_{D'} K_{(\psi,D-D')} \cdot
a^{D-D'}+\Lambda(a),
\end{eqnarray}
and $\displaystyle |\Lambda(a)|\leq \frac {C(n)} {a^n}$. Finally,
using $f(\lambda)= f^*(0) \big (|\lambda|^{-D}+C_{D'}
|\lambda|^{D'-D}\big ) +o\big ( |\lambda|^{D'-D} \big )$ when
$\lambda \to 0$, we obtain
\begin{eqnarray*}
\int_{-\sqrt{a}}^{\sqrt{a}}\Big ( f\big ( \frac u a \big ) - f^*(0)
\big (\Big | \frac u a \Big |^{-D} +C_{D'}\Big | \frac u a \Big
|^{D'-D}  \big )\Big )  |\widehat \psi (u)|^2\,
du&\\
&&\hspace{-7cm} =\int_{-\sqrt{a}}^{\sqrt{a}}\Big | \frac u a \Big
|^{D-D'} \Big (  f\big ( \frac u a \big ) - f^*(0) \big (\Big |
\frac u a \Big |^{-D} +C_{D'}\Big | \frac u a \Big |^{D'-D}  \big )
\Big ) |\widehat \psi (u)|^2 \Big |
\frac u a \Big |^{D'-D}\, du\\
&&\hspace{-7cm} =a^{D-D'} \int_{-\sqrt{a}}^{\sqrt{a}}g(u,a)
|\widehat \psi (u)|^2 |u|^{D'-D} du,
\end{eqnarray*}
with for all $u \in [ -\sqrt{a},\sqrt{a}]$, $g(u,a)\to 0$ when $a\to
\infty$. Therefore, from Lebesgue Theorem (checked from the
asymptotic behavior of $\widehat \psi$),
\begin{eqnarray}\label{som5}
\lim_{a\to \infty} a^{D-D'}\int_{-\sqrt{a}}^{\sqrt{a}}\Big ( f\big
( \frac u a \big ) - f^*(0) \big (\Big | \frac u a \Big |^{-D}
+C_{D'}\Big | \frac u a \Big |^{D'-D}  \big )\Big )  |\widehat
\psi (u)|^2\, du = 0.
\end{eqnarray}
As a consequence, from (\ref{som1}), (\ref{som2}), (\ref{som3}),
(\ref{som4}) and (\ref{som5}), the corollary is proven. $\Box$ \\
\dem[Proposition \ref{tlclog}] This proof can be decomposed into
three steps :{\bf Step 1}, {\bf Step 2} and {\bf Step 3}. \\
~\\
{\bf Step 1}. In this part, $\displaystyle{\frac N {a_N} \cdot \Cov
(\tilde T_N(r_ia_N), \tilde T_N(r_ja_N))_{1 \leq i,j \leq \ell}}$ is
proven to converge at an asymptotic  covariance matrix $\Gamma$.
First, for all $(i,j) \in \{1,\ldots,\ell\}^2$,
\begin{eqnarray}\label{covSN}
\Cov (\tilde T_N(r_ia_N), \tilde T_N(r_ja_N))=2 \frac 1
{[N/r_ia_N]}\frac 1 {[N/r_ja_N]} \sum_{p=1}^{[N/r_ia_N]}
\sum_{q=1}^{[N/r_ja_N]} \Big (\Cov(\tilde e(r_ia_N,p),\tilde
e(r_ja_N,q) \Big ) ^2,
\end{eqnarray}
because $X$ is a Gaussian process. Therefore, by considering only
$i=j$ and $p=q$, for $N$ and $a_N$ large enough,
\begin{eqnarray}\label{covSN1}
\Cov (\tilde T_N(r_ia_N), \tilde T_N(r_ia_N)) \geq \frac {1} {r_i}
\frac N {a_N}.
\end{eqnarray}
Now, for $(p,q) \in \{1,\ldots,[N/r_ia_N]\}\times
\{1,\ldots,[N/r_ia_N]\}$,
\begin{eqnarray*}
\Cov \big (\tilde e(r_ia_N,p),\tilde e(r_ja_N,q) \big )
\hspace{-3mm} &=& \hspace{-3mm}\frac {a_N^{1-D}(r_ir_j)^{(1-D)/2}}
{f^*(0) K_{(\psi,D)}}\frac 1 {r_i a_N} \frac 1 {r_j a_N}
\sum_{k=1}^{r_ia_N}\sum_{k'=1}^{r_j a_N} \psi\big (\frac k {r_ia_N}
\big ) \psi\big (\frac {k'} {r_ja_N} \big ) r\big
(k-k'+a_N(r_ip-r_jq)\big )\\
\hspace{-3mm} &=& \hspace{-3mm} \frac {a_N^{1-D}(r_ir_j)^{(1-D)/2}}
{f^*(0) K_{(\psi,D)}}\frac 1 {r_i a_N} \frac 1 {r_j a_N}
\sum_{k=1}^{r_ia_N}\sum_{k'=1}^{r_j a_N} \psi\big (\frac k {r_ia_N}
\big ) \psi\big (\frac {k'} {r_ja_N} \big ) \int_{-\pi}^\pi
\hspace{-3mm} d\lambda \,
f(\lambda)e^{-i\lambda(k-k'+a_N(r_ip-r_jq))}\\
\hspace{-3mm} &=& \hspace{-3mm} \frac {(r_ir_j)^{(1-D)/2}}
{a_N^{D}f^*(0) K_{(\psi,D)}}\frac 1 {r_i a_N} \frac 1 {r_j a_N}
\sum_{k=1}^{r_ia_N}\sum_{k'=1}^{r_j a_N} \psi\big (\frac k {r_ia_N}
\big ) \psi\big (\frac {k'} {r_ja_N} \big ) \int_{-\pi a_N}^{\pi
a_N}\hspace{-3mm} du \,  f\big (\frac u{a_N} \big )e^{-iu(\frac k
{a_N} -\frac {k'} {a_N}+r_ip-r_jq)}.
\end{eqnarray*}
Using the same expansion as in (\ref{Fourier2}), under Assumption
$W(\infty)$ the previous equality becomes, for all $n \in \N^*$,
\begin{eqnarray}
\nonumber  \left  | \Cov \big (\tilde e(r_ia_N,p),\tilde e(r_ja_N,q)
\big )- \frac {(r_ir_j)^{(1-D)/2}} {a_N^{D}f^*(0) K_{(\psi,D)}}
\int_{-\pi a_N}^{\pi a_N} \hspace{-3mm} du \, \widehat \psi
(ur_i)\overline{\widehat \psi} (ur_j) f\big (\frac u{a_N} \big
)e^{-iu(r_ip-r_jq)} \right |&& \\
\nonumber&& \hspace{-6cm} \leq \frac {C(n)} {a_N^{n+D}} \int_{-\pi
a_N}^{\pi a_N} \hspace{-3mm} du \,\big | \widehat \psi
(ur_i)\overline{\widehat \psi} (ur_j) f\big (\frac u{a_N} \big )
\big | \\
\nonumber && \hspace{-6cm} \leq \frac {C'(n)} {a_N^{n}}
\int_{-\infty}^{\infty } \hspace{-3mm} du \,|u|^{-D} \big | \widehat
\psi (ur_i)\overline{\widehat \psi} (ur_j) \big |\\
\label{ineq1} && \hspace{-6cm} \leq \frac {C''(n)} {a_N^{n}},
\end{eqnarray}
with $C(n),C'(n),C''(n)>0$ not depending on $a_N$ and due the
asymptotic behaviors of $\widehat \psi (u) $ when $u\to 0$ and $u\to
\infty$. Now, under Assumption A1,
\begin{eqnarray}
\nonumber&& \left | \int_{-\pi a_N}^{\pi a_N} \hspace{-3mm} du \,
\widehat \psi (ur_i)\overline{\widehat \psi} (ur_j) f\big (\frac
u{a_N} \big )e^{-iu(r_ip-r_jq)}-a_N f^*(0)\int_{-\pi}^{\pi }
\hspace{-3mm} du \, \frac {\widehat \psi (ur_ia_N)\overline{\widehat
\psi}
(ur_ja_N)}{|u|^D} e^{-iua_N(r_ip-r_jq)}\right | \\
\nonumber &&\hspace{8cm}  \leq a_N^{D-D'} f^*(0)C_{D'} \int_{-\pi
a_N}^{\pi a_N } \hspace{-3mm} du \,  \frac {\big | \widehat \psi
(ur_i)\overline{\widehat \psi} (ur_j) \big |}{|u|^{D-D'}}\\
\label{ineq2} &&\hspace{8cm}  \leq a_N^{D-D'} f^*(0)C_{D'}
\int_{-\infty}^{\infty } \hspace{-3mm} du \,  \frac {\big | \widehat
\psi (ur_i)\overline{\widehat \psi} (ur_j) \big |}{|u|^{D-D'}},
\end{eqnarray}
since  $\displaystyle \int_{-\infty}^{\infty } \hspace{-3mm} du \,
\frac {\big | \widehat \psi (ur_i)\overline{\widehat \psi} (ur_j)
\big |}{|u|^{D-D'}}<\infty$ from Assumption $W(\infty)$. Finally,
from (\ref{ineq1}) and (\ref{ineq2}), we have $C>0$ not depending on
$N$ such that for all $a_N\in \N^*$,
\begin{eqnarray}
\label{ineq3}\left  | \Cov \big (\tilde e(r_ia_N,p),\tilde
e(r_ja_N,q) \big )- \frac {a_N^{1-D} (r_ir_j)^{(1-D)/2}}
{K_{(\psi,D)}} \int_{-\pi}^{\pi } \hspace{-3mm} du \, \frac
{\widehat \psi (ur_ia_N)\overline{\widehat \psi} (ur_ja_N)}{|u|^D}
e^{-iua_N(r_ip-r_jq)} \right | \leq C \, a_N^{-D'}.
\end{eqnarray}
It remains to evaluate $\displaystyle a_N^{1-D} \int_{-\pi}^{\pi }
\hspace{-3mm} du \, \frac {\widehat \psi (ur_ia_N)\overline{\widehat
\psi} (ur_ja_N)}{|u|^D} e^{-iua_N(r_ip-r_jq)}=\int_{-\pi a_N}^{\pi
a_N} \hspace{-3mm} du \, \frac {\widehat \psi
(ur_i)\overline{\widehat \psi} (ur_j)}{|u|^D} e^{-iu(r_ip-r_jq)}$.
Thus, if $|r_ip-r_jq|\geq 1$, using an integration by parts,
\begin{eqnarray}
\nonumber \left | \int_{-\pi a_N}^{\pi a_N} \hspace{-3mm} du \,
\frac {\widehat \psi (ur_i)\overline{\widehat \psi} (ur_j)}{|u|^D}
e^{-iu(r_ip-r_jq)} \right | &=& \left | \frac 1 {-i(r_ip-r_jq)}\left
[\frac {\widehat \psi (ur_i)\overline{\widehat \psi}
(ur_j)}{u^D}e^{-iu(r_ip-r_jq)} \right]_{-\pi a_N}^{\pi  a_N} \right .   \\
\nonumber &+&\left . \frac 1 {i(r_ip-r_jq)}\int _{-\pi a_N}^{\pi
a_N} \hspace{-3mm} du \,\frac {\partial}{\partial u} \Big ( \frac
{\widehat \psi (ur_i)\overline{\widehat \psi}
(ur_j)}{u^D}\Big ) e^{-iu(r_ip-r_jq)}\right | \\
\nonumber  && \hspace{-6cm} \leq  \frac {1} {|r_ip-r_jq| }
\int_{-\infty}^{\infty}  \left (\frac D {|u|^{D+1}}\big |\widehat
\psi (ur_i)\overline{\widehat \psi} (ur_j)\big | + \frac 1
{|u|^{D}}\Big |  \frac {\partial}{\partial u} \Big (\widehat
\psi (ur_i)\overline{\widehat \psi} (ur_j) \Big )\Big | \right ) du \\
\label{equi1} && \hspace{-6cm} \leq C \frac {1} {|r_ip-r_jq| }
\end{eqnarray}
with $C<\infty$ not depending on $N$, since:
\begin{itemize}
\item $\widehat \psi (\pi
r_ia_N)\overline{\widehat \psi} (\pi r_ja_N)=\widehat \psi (-\pi
r_ia_N)\overline{\widehat \psi} (-\pi r_ja_N)$ and $ \sin(\pi
a_N(r_ip-r_jq))=0$;
\item from Assumption $W(\infty)$, $\limsup_{u\to 0} u^{-1}\,| \widehat
\psi (u)|<\infty$, $\limsup_{u\to 0} \big |\frac {\partial}{\partial
u} \widehat \psi (u)\big |<\infty$
$$
\Longrightarrow \limsup_{u\to 0} u^{-1}\,\Big | \frac
{\partial}{\partial u} \Big (\widehat \psi (ur_i)\overline{\widehat
\psi} (ur_j) \Big )\Big | <\infty;
$$
\item from Assumption $W(\infty)$, for all $n \in \N$,
$\sup_{u\in \R} (1+|u|)^n \,| \widehat \psi (u)|<\infty$ and
$\sup_{u\in \R} (1+|u|)^n \,\big |\frac {\partial}{\partial u}
\widehat \psi (u)\big |<\infty$.
\end{itemize}
Moreover, if $|r_ip-r_jq|= 0$, from Cauchy-Schwartz Inequality and
Property \ref{vard}, for $a_N$ large enough
\begin{eqnarray} \label{ineg_cov2}
\Big |\Cov \big (\tilde e(r_ia_N,p),\tilde e(r_ja_N,q) \big )\Big |
\leq  \Big (\E (\tilde e^2(r_ia_N,p)) \cdot \E (\tilde d^2(r_ja_N,q
)) \Big )^{1/2} \leq 2.
\end{eqnarray}
Therefore, using (\ref{ineq3}), (\ref{equi1}) and (\ref{ineg_cov2})
and the inequality $(x+y)^2 \leq 2 (x^2+y^2)$ for all $(x,y)\in
\R^2$, we have
 $C>0$ such that for $a_N$ large enough,
\begin{eqnarray}\label{ineq4}
\Cov^2 \big (\tilde e(r_ia_N,p),\tilde e(r_ja_N,q) \big ) \leq C
\Big ( \frac {1} {(1+|r_ip-r_jq|)^2}+ \frac {1} {a_N^{2D'}} \Big )
\end{eqnarray}
Hence, with (\ref{covSN}),
\begin{eqnarray*}
\Big | \Cov (\tilde T_N(r_ia_N), \tilde T_N(r_ja_N))\Big | \leq C
\frac 1 {[N/r_ia_N]}\frac 1 {[N/r_ja_N]} \sum_{p=1}^{[N/r_ia_N]}
\sum_{q=1}^{[N/r_ja_N]}\Big ( \frac {1} {(1+|r_ip-r_jq|)^2}+ \frac
{1} {a_N^{2D'}} \Big )
\end{eqnarray*}
But, from the theorem of comparison between sums and integrals,
\begin{eqnarray*}
\sum_{p=1}^{[N/r_ia_N]} \sum_{q=1}^{[N/r_ja_N]}
(1+|r_ip-r_jq|)^{-2}& \leq  & \frac 1 {r_ir_j} \int _0^{N/a_N}
\int _0^{N/a_N}\frac {du\, dv}{(1+|u-v|)^{2}} \\
& \leq  &\frac 2 {r_ir_j} \int _0^{N/a_N}  \frac {N/a_N \,
dw}{(1+w)^{2}}\\
& \leq & \frac 2 {r_ir_j} \cdot \frac N {a_N}.
\end{eqnarray*}
As a consequence, if $a_N$ is such that $\displaystyle \limsup
_{N\to \infty} \frac N {a_N} \frac 1 {a_N^{2D'}}<\infty$ then
$\displaystyle \limsup _{N\to \infty} \frac N {a_N} \Big | \Cov
(\tilde T_N(r_ia_N), \tilde T_N(r_ja_N))\Big |<\infty$. More
precisely, since this covariance is a sum of positive terms, if
$\displaystyle \limsup _{N\to \infty} \frac N {a_N} \frac 1
{a_N^{2D'}}=0$,
\begin{eqnarray}\label{covfinale}
\lim_{N\to \infty}\frac N {a_N}  \Big (\Cov (\tilde S_N(r_ia_N),
\tilde S_N(r_ja_N)) \Big) _{1\leq i,j \leq \ell}
=\Gamma(r_1,\cdots,r_\ell,\psi,D),
\end{eqnarray}
a non null (from (\ref{covSN1})) symmetric matrix with
$\Gamma(r_1,\cdots,r_\ell,\psi,D)=(\gamma_{ij})_{1\leq i,j\leq \ell}$
that can be specified. Indeed, from the previous computations, if
$\displaystyle \limsup _{N\to \infty} \frac N {a_N} \frac 1
{a_N^{2D'}}=0$,
\begin{eqnarray*}
\gamma_{ij}&=&\lim_{N\to \infty}\frac {8r_ir_j a_N} {N}
\sum_{p=1}^{[N/r_ia_N]} \sum_{q=1}^{[N/r_ja_N]}\Big ( \frac
{(r_ir_j)^{(1-D)/2}} { K_{(\psi,D)}} \int _0 ^\infty du \, \frac
{\widehat \psi (ur_i)\overline{\widehat \psi} (ur_j) }{u^D}\cos
(u(r_ip-r_jq)) \Big) ^2
\\
&=&\lim_{N\to \infty}\frac {8(r_ir_j)^{2-D} a_N} {K^2_{(\psi,D)}N}
\sum_{m=-[N/d_{ij}a_N]+1}^{[N/d_{ij}a_N]-1}(\frac N
{d_{ij}a_N}-|m|\big )\Big (  \int _0 ^\infty du \, \frac {\widehat
\psi (ur_i)\overline{\widehat \psi} (ur_j) }{u^D}\cos (u \,d_{ij} m)
\Big) ^2\\
&=&\frac {8(r_ir_j)^{2-D}} {K^2_{(\psi,D)}d_{ij}}
\sum_{m=-\infty}^{\infty}\Big ( \int _0 ^\infty \frac {\widehat \psi
(ur_i)\overline{\widehat \psi} (ur_j) }{u^D}\cos (u \,d_{ij} m)\, du
\Big) ^2,
\end{eqnarray*}
with $d_{ij}=GCD(r_i\,;\,r_j)$. Therefore, the matrix $\Gamma$
depends only on
on $r_1,\cdots,r_\ell,\psi,D$.\\
~\\
{\bf Step 2}.Generaly  speaking, the above result is not sufficient
to obtain the central limit theorem,
\begin{eqnarray}\label{TLC1}
\sqrt{\frac N {a_N}}  \Big (\tilde T_N(r_ia_N)-\E(\tilde
e^2(r_ia_N,0)\Big) _{1\leq i \leq \ell} \limiteloiN {\cal N}_\ell
(0,\Gamma(r_1,\cdots,r_\ell,\psi,D)).
\end{eqnarray}
However, each $\tilde T_N(r_ia_N)$ is a quadratic form of a Gaussian
process. {\it Mutatis mutandis}, it is exactly the same framework
({\it i.e.} a Lindeberg central limit theorem) as that of
Proposition 2.1 in Bardet (2000), and (\ref{TLC1}) is checked.
Moreover, if $(a_n)_n$ is such that $\displaystyle \limsup _{N\to
\infty} \frac N {a_N^{1+2D'} } =0$ then using the asymptotic
behavior of $\E(\tilde e^2(r_ia_N,0)$ provided in Property
\ref{vard},
$$
\sqrt{\frac N {a_N}}\Big (\E(\tilde e^2(r_ia_N,0)\Big ) \limitN 0.
$$
As a consequence, under those assumptions,
\begin{eqnarray}\label{TLC2}
\sqrt{\frac N {a_N}}  \Big (\tilde T_N(r_ia_N)-1\Big) _{1\leq i \leq
\ell} \limiteloiN {\cal N}_\ell
(0,\Gamma(r_1,\cdots,r_\ell,\psi,D)).
\end{eqnarray}
~\\
{\bf Step 3}. The logarithm function $(x_1,..,x_\ell) \in
(0,+\infty)^\ell \mapsto (\log x_1,..,\log x_m)$ is ${\cal C}^2$
on $(0,+\infty)^\ell$. As a consequence, using the Delta-method,
the central limit theorem (\ref{CLTSN}) for the vector
$\displaystyle{\left (\log {\tilde T}_N(r_ia_N) \right ) _{1 \leq
i \leq \ell}}$ follows with the same asymptotical covariance
matrix $\Gamma(r_1,\cdots,r_\ell,\psi,D)$ (because the Jacobian
matrix of the function in $(1,..,1)$ is the identity matrix). $\Box$\\
~\\
\dem[Proposition \ref{tlclog2}] There is a perfect identity between
this proof and that of Proposition \ref{tlclog}, both of which are
based on the   approximations of Fourier transforms provided in the proof of Property \ref{vard'}. $\Box$\\
\dem[Corollary \ref{cor3}] It is clear that $X'_t=X_t+P_m(t)$ for
all $t\in \Z$, with $X=(X_t)_t$ satisfying Proposition
\ref{tlclog} and \ref{tlclog2}. But, any wavelet coefficient of
$(P_m(t))_t$ is obviously null from
the assumption on $\psi$. Therefore the statistic $\widehat T_N$ is the same for $X$ and $X'$. $\Box$\\
\dem[Proposition \ref{hatalpha}] Let $\varepsilon>0$ be a fixed
positive real number, such that $\alpha^*+\varepsilon<1$. \\
~\\
{\bf I.} First, a bound of $\Pr(\widehat \alpha_N \leq
\alpha^*+\varepsilon)$ is provided. Indeed,
\begin{eqnarray}
\nonumber \Pr\big (\widehat \alpha_N \leq \alpha^*+\varepsilon \big
) & \geq & \Pr\Big (\widehat Q_N(\alpha^*+\varepsilon/2) \leq
\min_{\alpha \geq \alpha^*+\varepsilon~\mbox{and}~\alpha \in {\cal
A}_N}\widehat Q_N(\alpha)\Big ) \\
\nonumber & \geq & 1 - \Pr\Big (\bigcup_{\alpha \geq
\alpha^*+\varepsilon~\mbox{and}~\alpha \in {\cal A}_N}\widehat
Q_N(\alpha^*+\varepsilon/2) > \widehat Q_N(\alpha)\Big )  \\
 \label{QN} & \geq &1 -  \sum_{k =[(\alpha^*+\varepsilon)\log N]}^{\log
[N/\ell]}\Pr\Big (\widehat Q_N(\alpha^*+\varepsilon/2) > \widehat
Q_N\big (\frac {k}{\log N} \big)\Big ).
\end{eqnarray}
But, for $\alpha\geq \alpha^*+1$,
\begin{eqnarray*}
\Pr\Big (\widehat Q_N(\alpha^*+\varepsilon/2) > \widehat
Q_N(\alpha)\Big )= \Pr\Big (\Big \| P_N(\alpha^*+\varepsilon/2)
\cdot  Y_N(\alpha^*+\varepsilon/2) \Big \|^2  > \Big \|
P_N(\alpha) \cdot Y_N(\alpha)\Big \|^2 \Big )
\end{eqnarray*}
with $P_N(\alpha)=I_\ell -A_N(\alpha) \cdot  \big (
A_N'(\alpha)\cdot A_N(\alpha) \big )^{-1} \cdot A_N(\alpha)$ for all
$\alpha \in (0,1)$, {\it i.e.} $P_N(\alpha)$ is the matrix of an
orthogonal projection on the orthogonal subspace (in $\R^\ell$)
generated by $A_N(\alpha)$ (and $I_\ell$ is the identity matrix in
$\R^\ell$). From the expression of $A_N(\alpha)$, it is obvious that
for all $\alpha \in (0,1)$,
$$
P_N(\alpha)=P=I_\ell -A \cdot  \big ( A'\cdot A \big )^{-1} \cdot
A,
$$
with the matrix $\displaystyle{A= \left
( \begin{array}{cc}\log (r_1) & 1  \\
: & : \\ \log (r_\ell ) & 1  \end{array}\right )}$ as  in
Proposition \ref{tlcd}. Thereby,
\begin{eqnarray*}
\Pr\Big (\widehat Q_N(\alpha^*+\varepsilon/2) > \widehat
Q_N(\alpha)\Big )&= & \Pr\Big (\Big \| P \cdot
Y_N(\alpha^*+\varepsilon/2) \Big \|^2  > \Big \| P \cdot
Y_N(\alpha)\Big \|^2 \Big ) \\
&= & \Pr\left ( {\displaystyle  \Big \| P \cdot\sqrt{\frac N
{N^{\alpha^*+\varepsilon/2}}} \,Y_N(\alpha^*+\varepsilon/2) \Big
\|^2}
> N^{\alpha-(\alpha^*+\varepsilon/2)} {\displaystyle \Big \| P \cdot\sqrt{\frac N {N^{\alpha}}} \,
Y_N(\alpha) \Big \|^2}\right )\\
& \leq  &  \Pr\left ( V_N(\alpha^*+\varepsilon/2)
> N^{(\alpha-(\alpha^*+\varepsilon/2))/2} \right) +  \Pr\left ( V_N(\alpha)
\leq  N^{-(\alpha-(\alpha^*+\varepsilon/2))/2} \right)
\end{eqnarray*}
with $\displaystyle{V_N(\alpha)=  \Big \| P \cdot\sqrt{\frac N
{N^{\alpha}}} \,Y_N(\alpha) \Big \|^2}$ for all $\alpha \in (0,1)$.
From Proposition \ref{tlclog}, for all $\alpha >\alpha^*$, the
asymptotic law of $\displaystyle P \cdot\sqrt{\frac N {N^{\alpha}}}
\, Y_N(\alpha)$ is a Gaussian law with covariance matrix $P\cdot
\Gamma \cdot P'$. Moreover, the rank of the matrix is $P\cdot \Gamma
\cdot P'$ is $\ell-2$ (this is the rank of $P$) and we have

$0<\lambda_-$, not depending on $N$) such that $P\cdot \Gamma \cdot
P'-\lambda_- P \cdot P' $ is a non-negative matrix ($0<\lambda_- <
\min\{ \lambda \in \mbox{Sp}(\Gamma)\}$). As a consequence, for a
 large enough $N$,
\begin{eqnarray*}
\Pr\left ( V_N(\alpha) \leq
N^{-(\alpha-(\alpha^*+\varepsilon/2))/2} \right) &\leq& 2 \cdot
\Pr\left ( V_- \leq N^{-(\alpha-(\alpha^*+\varepsilon/2))/2}
\right)
\\
& \leq & \frac 1 {2^{\ell/2-2}\Gamma(\ell/2)}\cdot\Big (  \frac N
{\lambda_-} \Big ) ^{-(\frac \ell 2 -1)\frac {
(\alpha-(\alpha^*+\varepsilon/2))} 2},
\end{eqnarray*}
with $V_- \sim \lambda_- \cdot \chi^2(\ell-2)$. Moreover, from
Markov inequality,
\begin{eqnarray*}
\Pr\left ( V_N(\alpha^*+\varepsilon/2)
> N^{(\alpha-(\alpha^*+\varepsilon/2))/2} \right) & \leq & 2 \cdot \Pr\left (
\exp(
\sqrt{V_+}) > \exp \big ( N^{(\alpha-(\alpha^*+\varepsilon/2))/4}\big ) \right)\\
& \leq & 2 \cdot \E (\exp( \sqrt{V_+})) \cdot \exp \big (-
N^{(\alpha-(\alpha^*+\varepsilon/2))/4}\big )
\end{eqnarray*}
with $V_+ \sim \lambda_+ \cdot \chi^2(\ell-2)$ and $\lambda_+ >
\max\{ \lambda \in \mbox{Sp}(\Gamma)\}>0$. Like $\E (\exp(
\sqrt{V_+}))<\infty$ does not depend on $N$, we obtain that $M_1>0$
not depending on $N$, such that for  large enough $N$,
\begin{eqnarray*}
\Pr\Big (\widehat Q_N(\alpha^*+\varepsilon/2) > \widehat
Q_N(\alpha)\Big )  \leq M_1 \cdot N ^{-(\frac \ell 2 -1)\frac {
(\alpha-(\alpha^*+\varepsilon/2))} 2},
\end{eqnarray*}
and therefore, the inequality (\ref{QN}) becomes, for $N$ large
enough,
\begin{eqnarray}
\nonumber \Pr\big (\widehat \alpha_N \leq \alpha^*+\varepsilon \big
) & \geq &1 - M_1\cdot \sum_{k =[(\alpha^*+\varepsilon)\log N]}^{\log
[N/\ell]} N ^{-\frac {(\ell-2)} 4\Big (\big (\frac {k}{\log N}\big
) -(\alpha^*+\varepsilon/2)\Big ) }\\
\label{borne1} & \geq &1 - M_1 \cdot \log N \cdot N^{-\frac {(\ell-2)}
{12} \varepsilon }.
\end{eqnarray}
{\bf II.} Secondly, a bound of $\Pr(\widehat \alpha_N \geq
\alpha^*-\varepsilon)$ is provided. Following the above arguments
and notations ,
\begin{eqnarray}
\nonumber \Pr\big (\widehat \alpha_N \geq \alpha^*-\varepsilon \big
) & \geq & \Pr\Big (\widehat Q_N(\alpha^*+ \frac
{1-\alpha^*}{2\alpha^*} \varepsilon) \leq \min_{\alpha \leq
\alpha^*-\varepsilon~\mbox{and}~\alpha \in {\cal
A}_N}\widehat Q_N(\alpha)\Big ) \\
\label{QN2} & \geq &1 -  \sum_{k =2}^{[(\alpha^*-\varepsilon) \log
N]+1}\Pr\Big (\widehat Q_N(\alpha^*+\frac {1-\alpha^*} {2 \alpha^*}
\varepsilon) > \widehat Q_N\big (\frac {k}{\log N} \big)\Big ),
\end{eqnarray}
and as above,
\begin{eqnarray}
\nonumber  \Pr\Big (\widehat Q_N(\alpha^*+\frac {1-\alpha^*} {2
\alpha^*} \varepsilon) > \widehat Q_N(\alpha)\Big )&& \\
\label{TLC7} && \hspace{-4.5cm} =  \Pr\left ( {\displaystyle  \Big
\| P \cdot\sqrt{\frac N {N^{\alpha^*+\frac {1-\alpha^*} {2 \alpha^*}
\varepsilon}}} \,Y_N(\alpha^*+\frac {1-\alpha^*} {2 \alpha^*}
\varepsilon) \Big \|^2}
> N^{\alpha-(\alpha^*+\frac {1-\alpha^*} {2
\alpha^*} \varepsilon)} {\displaystyle \Big \| P \cdot\sqrt{\frac N
{N^{\alpha}}} \, Y_N(\alpha) \Big \|^2}\right ).
\end{eqnarray}
Now, in the case $a_N=N^\alpha$ with $\alpha \leq \alpha^*$, the
sample variance of wavelet coefficients is biased. In this case,
from the relation of Corollary \ref{cor1} under Assumption A1',
\begin{eqnarray*}
\Big ( Y_N(\alpha) \Big )_{1\leq i \leq \ell}=\Big ( \frac { C_{D'}
K_{(\psi,D-D'))}} {f^*(0)K_{(\psi,D)}} (iN^{\alpha})^{-
D'}(1+o_i(1)) \Big )_{1\leq i \leq \ell} + \Big ( \sqrt{\frac
{N^{\alpha}}N}\cdot \varepsilon_N(\alpha) \Big )_{1\leq i \leq
\ell},
\end{eqnarray*}
with $o_i(1) \to 0$ when $N\to \infty $ for all $i$ and $\E
(Z_N(\alpha))=0$. As a consequence, for  large enough $N$,
\begin{eqnarray*}
\Big \| P \cdot\sqrt{\frac N {N^{\alpha}}} \, Y_N(\alpha) \Big
\|^2&=&\Big \| P \cdot \varepsilon_N(\alpha) \Big
\|^2+N^{\frac{\alpha^*-\alpha}{\alpha^*}}\Big \| P \cdot \Big (
\frac { C_{D'} K_{(\psi,D-D'))}} {f^*(0)K_{(\psi,D)}}
i^{- D'}(1+o_i(1)) \Big )_{1\leq i \leq \ell} \Big \|^2 \\
& \geq &  D \cdot N^{\frac{\alpha^*-\alpha}{\alpha^*}},
\end{eqnarray*}
with $D >0$,  because the vector $(i^{-D'})_{1\leq i \leq \ell}$
is not in the orthogonal subspace of the subspace generated by the
matrix $A$. Then, the relation (\ref{TLC7}) becomes,
\begin{eqnarray*}
\Pr\Big (\widehat Q_N(\alpha^*+\frac {1-\alpha^*} {2 \alpha^*}
\varepsilon) > \widehat Q_N(\alpha)\Big ) &\leq & \Pr\Big (
{\displaystyle  \Big \| P \cdot\sqrt{\frac N {N^{\alpha^*+\frac
{1-\alpha^*} {2 \alpha^*} \varepsilon}}} \,Y_N(\alpha^*+\frac
{1-\alpha^*} {2 \alpha^*} \varepsilon) \Big \|^2} \geq D \cdot
N^{\alpha-(\alpha^*+\frac {1-\alpha^*} {2
\alpha^*} \varepsilon)} \cdot N^{\frac {\alpha^*-\alpha}{\alpha^*}}\Big ) \\
&\leq & \Pr\Big ( V_+ \geq D \cdot N^{\frac {1-\alpha^*} {2
\alpha^*}(2(\alpha^*-\alpha)- \varepsilon)}\Big
) \\
& \leq & M_2 \cdot N^{-(\frac \ell 2 -1)\frac {1-\alpha^*} {2
\alpha^*}\varepsilon},
\end{eqnarray*}
with $M_2 >0$, because $V_+ \sim \lambda_+ \cdot \chi^2(\ell-2)$ and
$\displaystyle{\frac {1-\alpha^*} {2 \alpha^*}(2(\alpha^*-\alpha)-
\varepsilon)\geq \frac {1-\alpha^*} {2 \alpha^*}\varepsilon}$ for
all $\alpha \leq \alpha^*-\varepsilon$. Hence, from the inequality
(\ref{QN2}), for large enough $N$ ,
\begin{eqnarray}
\label{borne2} \Pr\big (\widehat \alpha_N \geq \alpha^*-\varepsilon
\big ) \geq 1 -M_2 \cdot \log N \cdot N^{-(\frac \ell 2 -1)\frac
{1-\alpha^*} {2 \alpha^*}\varepsilon}.
\end{eqnarray}
The inequalities (\ref{borne1}) and (\ref{borne2}) imply that
$\displaystyle{~\Pr\big (|\widehat \alpha_N - \alpha| \geq
\varepsilon \big ) \limitN 0}$. $\Box$ \\
~\\
\dem[Theorem \ref{tildeD}] The central limit theorem of
(\ref{CLTD2}) can be established from the following arguments.
First, $\Pr (\tilde \alpha_N>\alpha^*) \limitN 1$. Following the
previous proof, there is for all $\varepsilon>0$,
\begin{eqnarray*}
\Pr\big (\widehat \alpha_N \geq \alpha^*-\varepsilon \big ) \geq 1
-M_2 \cdot \log N \cdot N^{-(\frac \ell 2 -1)\frac {1-\alpha^*} {2
\alpha^*}\varepsilon}.
\end{eqnarray*}
Consequently, if $\displaystyle{\varepsilon_N=\lambda \cdot \frac
{\log \log N}{\log N}}$ with $\displaystyle{\lambda> \frac 2
{(\ell-2)D'} }$ then,
\begin{eqnarray*}
\Pr\big (\widehat \alpha_N \geq \alpha^*-\varepsilon_N \big ) &\geq
&  1 -M_2 \cdot \log N \cdot N^{-\lambda \frac {(\ell- 2)D'} 2  \cdot
\frac {\log \log N}{\log N}} \\
&\geq & 1 -M_2 \cdot \big (\log N \big )^{1-\lambda \frac {(\ell-
2)D'} 2}\\
\Longrightarrow ~~ && \Pr\big (\widehat \alpha_N +\varepsilon_N
\geq \alpha^*\big ) \limitN  1.
\end{eqnarray*}
Now, from Corollary \ref{hatD'}, $\widehat {D'}_N \limiteprobaN
D'$. Therefore, $\displaystyle{\Pr\big (\widehat {D'}_N \leq \frac
4 3 D' \big ) \limitN 1}$. Thus, with $\displaystyle{\lambda \geq
\frac 9 {4(\ell-2)D'} }$,  $\displaystyle{\Pr\big (\tilde \alpha_N
+\big (\varepsilon_N -\frac 3 {(\ell-2)\widehat {D'}_N  }\cdot
\frac {\log \log N}{\log N}\big ) \geq \alpha^*\big ) \limitN 1}$
which implies $\Pr (\tilde
\alpha_N>\alpha^*) \limitN 1$. \\
Secondly, for $x \in _R$,
\begin{eqnarray*}
\lim_{N\to \infty} \Pr \Big (\sqrt{\frac N {N^{\tilde
\alpha_N}}}\big(\tilde D_N - D \big)\leq x  \Big )&=&\lim_{N\to
\infty} \Pr \Big (\sqrt{\frac N {N^{\tilde \alpha_N}}}\big(\tilde
D_N - D \big)\leq x \bigcap \tilde \alpha_N>\alpha^* \Big ) \\
&& \hspace{1cm}+ \lim_{N\to \infty} \Pr \Big (\sqrt{\frac N
{N^{\tilde \alpha_N}}}\big(\tilde D_N - D \big)\leq x \bigcap
\tilde \alpha_N\leq \alpha^* \Big )\\
&=&\lim_{N\to \infty} \int _{\alpha^*}^1 \Pr \Big (\sqrt{\frac N
{N^{\alpha}}}\big(\tilde
D_N - D \big)\leq x  \Big )f_{\widehat \alpha_N}(\alpha)\,d\alpha  \\
&=&\lim_{N\to \infty} \Pr \Big (Z_\Gamma \leq x  \Big ) \cdot \int
_{\alpha^*}^1
f_{\widehat \alpha_N}(\alpha)\,d\alpha  \\
& =& \Pr \Big (Z_\Gamma \leq x  \Big ),
\end{eqnarray*}
with $f_{\widehat \alpha_N}(\alpha)$ the probability density
function of $\widehat \alpha_N$ and $Z_\Gamma \sim \mathcal{N}(0\,
; \, (A' \cdot A ) ^{-1}\cdot A' \cdot \Gamma \cdot A \cdot (A'
\cdot A ) ^{-1})$. \\
\\
To prove the second part of (\ref{CLTD2}), we infer  deduces from
above that
$$
\Pr \Big (\alpha^*<\tilde \alpha_N <\alpha^*+\frac 3
{(\ell-2)\widehat {D'}_N  }\cdot \frac {\log \log N}{\log N}+ \mu
\cdot \frac {\log \log N}{\log N}\Big ) \limitN 1,
$$
with $\mu
>\frac {12}{\ell-2}$. Therefore,  $\nu <\frac 4
{(\ell-2)D'}+\frac {12}{\ell-2}$,
$$
\Pr\Big (N^{\alpha^*} <N^{\tilde \alpha_N} <N^{\alpha^*}\cdot
(\log N )^\nu \Big ) \limitN 1.
$$
This inequality and the previous central limit theorem result in :
 for all $\rho>\nu/2$, and $\varepsilon>0$,
\begin{eqnarray*}
\Pr \Big (\frac {N^{\frac {D'}{1+2D'}}}{(\log N)^{\rho}} \cdot
\big|\tilde D_N - D \big|
>\varepsilon \Big )&=&\Pr \Big (\frac {N^{\frac 1 2 (\widehat \alpha_N -\alpha^*)}}{(\log N)^{\rho}}
\cdot \sqrt{\frac N {N^{\tilde \alpha_N}}}\big|\tilde D_N - D
\big|
>\varepsilon  \Big ) \\
&\limitN& 0. ~\Box
\end{eqnarray*}
{\bf Acknowledgments.} The authors are very grateful to
anonymous referees for many relevant suggestions and corrections
that strongly improve the content of the paper.
\bibliographystyle{amsalpha}

\end{document}